%
%
%

\input amstex
\documentstyle{amsppt}
\PSAMSFonts

\pageheight{18.5cm}
\magnification=\magstep1
\frenchspacing


\loadbold

\def\phi{\varphi}

\def\nin{\newline\indent}

\def\CC{{\Bbb C}}
\def\FF{{\Bbb F}}
\def\GG{{\Bbb G}}
\def\HH{{\Bbb H}}
\def\QQ{{\Bbb Q}}

\def\NN{{\Bbb N}}

\def\ZZ{{\Bbb Z}}
\def\LL{{\Bbb L}}

\def\GL{\operatorname{GL}}
\def\SL{\operatorname{SL}}

\def\FFF{{\Cal F}}

\def\RRR{{\Cal R}}
\def\OOO{{\Cal O}}
\def\QQQ{{\Cal Q}}
\def\VVV{{\Cal V}}

\def\PPP{{\Cal P}}

\def\GGG{{\Cal G}}
\def\KKK{{\Cal K}}

\def\qqed{{\hfill\hfill\qed}}

\def\se#1,#2,#3;{(ab)^{#1}\, a_x^{\,#2}b_x^{\,#3}}

\def\see#1,#2,#3;{(ab)^{#3} (bc)^{#1} (ca)^{#2}\,
\ell_a^{p-#2-#3}\ell_b^{q-#1-#3}\ell_c^{r-#1-#2}}

\def\seep#1,#2,#3;{(ab)^{#3} (bc)^{#1} (ca)^{#2}\,
\ell_a^{p'-#2-#3}\ell_b^{q'-#1-#3}\ell_c^{r'-#1-#2}}

\def\seed#1,#2,#3;{(ab)^{#3} (bc)^{#1} (ca)^{#2}\,
a_x^{p_a-#2-#3}b_x^{p_b-#1-#3}c_x^{p_c-#1-#2}}

\def\nse#1,#2,#3,#4,#5,#6;{(#1 #2)(#3 #4){#5}_x{#6}_x}
\def\nsee#1,#2,#3,#4,#5,#6;{(#1 #2)(#3 #4)#5_x#6_x}
\def\nsse#1,#2,#3,#4,#5,#6,#7,#8;{(#1 #2)(#3 #4)(#5 #6)(#7 #8)}

\def\tv#1,#2,#3;{[#1,#2]_{#3}}

\TagsOnRight

\topmatter
\title Generic free resolutions and root systems
\endtitle

\author Jerzy Weyman \endauthor
\address Department of Mathematics, University of Connecticut
\nin Storrs, CT 06269, USA \endaddress
\email jerzy.weyman\@uconn.edu \endemail

\thanks The  author was partially supported by NSF grants DMS-0901185 and DMS 1400740. He also acknowledges the support of Alexander von Humboldt Foundation.
\endthanks

\rightheadtext{Resolutions and root systems}
\abstract In this paper I give an explicit construction of the generic rings ${\hat R}_{gen}$ for  
free resolutions of length 3 over Noetherian commutative 
 $\CC$-algebras.
The key role is played by the defect Lie algebra introduced in \cite{W89}.
The defect algebra turns out to be a parabolic subalgebra in a Kac-Moody Lie algebra associated
to the graph $T_{p,q,r}$ corresponding to the format of the resolution.
The  ring ${\hat R}_{gen}$ is Noetherian if and only if the graph $T_{p,q,r}$ corresponding to a given format 
 is a Dynkin diagram. In such case ${\hat R}_{gen}$ has rational singularities so it is Cohen-Macaulay.
The  ring ${\hat R}_{gen}$ is a deformation of a commutative ring ${\hat R}_{spec}$
which has a structure of a multiplicity free module
over a product of Kac-Moody Lie algebra corresponding to the graph $T_{p,q,r}$ and a product of two
general linear Lie algebras.
\endabstract
\endtopmatter

\document

\head  \S1. Introduction \endhead

In this paper I consider a problem of existence of  generic free resolutions of length three.

Let us look at  free acyclic complexes  $\FF_\bullet$ (i.e complexes whose only nonzero homology group is $H_0(\FF_\bullet)$) of the form
$$\FF_\bullet:\ 0\rightarrow F_3\rightarrow F_2\rightarrow F_1\rightarrow F_0$$
over  commutative Noetherian rings $R$,
with $rank\ F_i=f_i$ ($0\le i\le 3$), $rank (d_i)=r_i$ ($1\le i\le 3$). The quadruple $(f_0, f_1, f_2, f_3)$ is {\it the format} of the complex $\FF_\bullet$. 
We always have $f_i=r_i+r_{i+1}$ ($0\le i\le 3$), with the convention that $r_4=0$.

For the resolutions of such format $(f_0, f_1, f_2, f_3)$ we say that a pair $(R_{gen}, \FF^{gen}_\bullet)$ where $R_{gen}$ is a commutative ring and $\FF^{gen}_\bullet$ is an acyclic free complex of format $(f_0, f_1, f_2, f_3)$ over $R_{gen}$ is {\it a generic resolution} of this format
if two conditions are satisfied:

\item{a)} The complex $\FF^{gen}_\bullet$ is acyclic over $R_{gen}$,
\item{b)} For every acyclic free complex $\GG_\bullet$ over a Noetherian ring $S$ there exists a ring homomorphism $\phi:R_{gen}\rightarrow S$ such that
$$\GG_\bullet =\FF^{gen}_\bullet\otimes_{R_{gen}}S.$$

Of particular interest is whether the ring $R_{gen}$ is Noetherian, because it can be shown quite easily that a non-Noetherian (non-unique) generic pair always exists(\cite{Br84}).

For complexes of length $2$ the existence of the pairs $(R_{gen}, \FF^{gen}_\bullet)$ was established by Hochster (\cite{H75}). He also proved that this generic ring is Noetherian.
Later the generic rings for resolutions of length two were described very explicitly (see \cite{PW90} and references therein). This case is, however special, as only for the case of the resolution of length two the ring $R_{gen}$ is unique.

For free resolutions of length three the situation is much more difficult. We limit ourselves to considering the case of characteristic zero.
Thus all commutative rings considered in this paper are $\CC$-algebras, i.e. we work over an algebraically closed field of characteristic zero. 
The results are true also for $\QQ$-algebras as the
representations of the group $GL_n$ are defined over $\QQ$.

In \cite{W89} I constructed for each format of length three  certain pair $({\hat R}_{gen}, \hat{\FF}^{gen}_\bullet)$ obtained  by a specific process of killing cycles in a generic complex. I conjectured that the pair $({\hat R}_{gen}, \hat{\FF}^{gen}_\bullet)$ is  generic for each format of length three.

The key point in the construction was defining a family of certain  graded Lie algebra $\LL (r, E, F)$ (defect Lie algebras)
depending on a natural number $r$ and  two vector spaces $E, F$ over $\CC$. 
The defect Lie algebra related to the resolution $\FF_\bullet$ is the Lie algebra $\LL (r_1+1, F_3,F_1)$.

In \cite{W89} I proved that the assertion that the pair $({\hat R}_{gen}, \hat{\FF}^{gen}_\bullet)$ gives a generic pair is true  provided a certain family of three term complexes over the enveloping algebra of the corresponding defect Lie algebra $\LL (r, E, F)$ is exact at a middle term.

\vskip .3cm

In the present paper this program is completed. 

I associate to the format $(f_0, f_1, f_2, f_3)$  the graph $T_{p,q,r}$ with $p, q, r$ given by the rule $r_3=r-1$, $r_2=q+1$, $r_1=p-1$.

The key new observation is that the defect Lie algebra $\LL(p, E, F)$ coresponding to our format is the parabolic subalgebra associated to certain grading of the Kac-Moody Lie algebra  ${\goth g}(T_{p,q,r})$. This grading is determined by choosing a particular simple root of ${\goth g}(T_{p,q,r})$. The complexes constructed in \cite{W89} turn out to be parts of the parabolic BGG complexes associated to this grading on  ${\goth g}(T_{p,q,r})$. This shows immediately that these complexes are exact at prescribed places and thus it proves that the pair $({\hat R}_{gen}, \hat{\FF}^{gen}_\bullet)$ indeed is a generic pair  for every format of length three. 
 
 Moreover, the link with Kac-Moody Lie algebra allows to prove much more precise results.  It turns out that the generic ring ${\hat R}_{gen}$ has a deformation to
 a commutative ring ${\hat R}_{spec}$ which has a  structure of a module over a bigger Kac-Moody Lie algebra ${\goth g}(T_{p,q,r})\times {\goth gl}(F_2)\times{\goth gl}(F_0)$. This allows to prove several properties of ${\hat R}_{gen}$.
 
 We show that the ring ${\hat R}_{gen}$ is Noetherian whenever the Lie algebra $\LL$ is finite dimensional.
 Thus the  generic ring ${\hat R}_{gen}$ is Noetherian if and only if the diagram $T_{p,q,r}$ is a Dynkin diagram. Moreover, in Noetherian cases the ring ${\hat R}_{gen}$ has rational singularities and so it is Cohen-Macaulay.
 
 The paper is organized as follows. 
 
 In section 2 we recall the basic notions on free resolutions and on the grade of ideals in the non-Noetherian rings (true grade in the sense of Northcott) we will use throughout.
 
 Section 3  recalls the basic results on Kac-Moody Lie algebras.
 In section 4  more precise formulas for the Kac-Moody Lie algebras associated to the graphs $T_{p,q,r}$ are worked out. The reader is advised to skip these two sections in the first reading.
 
 Section 5 describes the geometric approach using higher derived image functors  to proving acyclicity of complexes.
 
 In section 6 we discuss the structure of the rings $R_a$ generated by the entries of maps in a generic complex and the Buchsbaum-Eisenbud multipliers.
 These rings are the staring point of our construction. We make explicit the representation structure of the rings $R_a$ and of the homology of the 
 generic complex $\FF_a$ over $R_a$.
 
 In section 7 we define the higher structure maps $p_i$. These are the cycles we kill to get the generic ring ${\hat R}_{gen}$.
 We also show the connection of these maps with the defect Lie algebra.
 
 Section 8 summarizes  the approach from \cite{W89}. We  recall there that genericity of ${\hat R}_{gen}$ would follow from exactness of some three term complexes. 
 
 In section 9  the main results are proved. The  identification of  the defect algebra with the parabolic algebra in the Kac-Moody Lie algebra ${\goth g}(T_{p,q,r})$ is carried out.
 The three terms complexes become parts of parabolic BGG resolutions of ${\goth g}(T_{p,q,r})$ so they are exact in needed places.
 
 In section 10 we analyze the representation structure, generators and relations of the rings ${\hat R}_{gen}$, and discuss its deformation to ${\hat R}_{spec}$. We also prove that in cases these rings are Noetherian they have rational singularities.

 In section 11 we give precise descriptions of the pairs $({\hat R}_{gen}, \hat{\FF}^{gen}_\bullet)$ in the simplest cases. It means we give a ''naive" description of the generators as entries in the matrices of some cycles whose existence one could deduce just from the universality property of ${\hat R}_{gen}$.
 
 Section 12 contains  some applications, including questions and conjectures about possible variations of the present approach.
 We are interested in particular in extending our results to perfect resolutions. The triality for the graphs $T_{p,q,r}$ suggests a possible connection with
 linkage theory.

\bigskip\bigskip

\head \S2. Background on finite free resolutions.\endhead

\bigskip

Throughout the paper we work with  complexes 
$$\FF_\bullet :0\rightarrow F_n{\buildrel{d_n}\over\rightarrow} F_{n-1}\rightarrow \ldots \rightarrow F_i{\buildrel{d_i}\over\rightarrow}F_{i-1}\rightarrow\ldots\rightarrow F_1{\buildrel{d_1}\over\rightarrow} F_0$$
of the free modules over a commutative Noetherian ring $R$. The map $d_i$ can be thought of  (after choosing bases in modules $F_i$) as a matrix with entries in $R$.
We define the rank $r_i$ of $d_i$ as the size of the biggest non-vanishing minor of that matrix. It is routine to check that this notion does not depend on the choice of bases in our free modules.
The ideal of $r_i\times r_i$ minors of $d_i$ is denoted $I(d_i)$.

We assume that $d_i :F_i\rightarrow F_{i-1}$ has rank $r_i$, and that
$$f_i :=rank\ F_i =r_i+r_{i+1}$$
with the convention $r_{n+1}=0, r_0\ge0$. We fix the numbers $f_i, r_i$.
We will say that a complex with these ranks has {\it a format} $(f_0, \ldots , f_n)$. We assume throughout that the format is 
valid, i.e. the numbers $r_i$ satisfying the above conditions exist.

All of this is motivated by the Buchsbaum-Eisenbud exactness criterion.

\proclaim{Theorem 2.1}(\cite{BE73})The complex $\FF$ over a Noetherian ring $R$ is acyclic if and only if 
\item{a)} $r_i+r_{i+1}= rank\ F_i$ for all $i=1,\ldots ,n$,
\item{b)} $depth_R (I(d_i)\ge i$ for $i=1,\ldots ,n$.
\endproclaim

A finite free resolution over $R$ is a free acyclic complex  $\FF$.

Buchsbaum and Eisenbud proved a fundamental structure theorem (the First Structure Theorem from \cite{BE74}).

\proclaim{Theorem 2.2}
Let $\FF$ be a finite free resolution over a Noetherian ring $R$. Then there exist unique maps $a_i:R\rightarrow\bigwedge^{r_i}F_i$
such that
\item{a)} $a_n=\bigwedge^{r_n}d_n$,
\item{b)} We have  factorizations
$$\matrix
\bigwedge^{r_i}F_i=\bigwedge^{r_{i+1}}F_i^*&&\buildrel{\bigwedge^{r_i}d_i}\over\rightarrow&&\bigwedge^{r_i}F_{i-1}\\
&\searrow\ a_{i+1}^*&&\nearrow\ a_i&\\
&&R\endmatrix$$
\endproclaim

The coordinates of maps $a_i$ are the Pl\"ucker coordinates of the images of maps $d_i$. In this formulation we identify the complementary exterior powers 
$\bigwedge^{r_i}F_i$ and $\bigwedge^{r_{i+1}}F_i^*$ so the diagram above is $\prod_iSL(F_i)$ equivariant. One can easily write down more precise $\prod_iGL(F_i)$ equivariant version, and we do it in section 6.
We refer to the maps $a_i$ as to Buchsbaum-Eisenbud multipliers.

Later Bruns generalized this result to prove the existence of the maps $a_i$ for all free complexes $\FF$ for which the rank condition $r_i+r_{i+1}=rank(F_i)$ for $i=1,\ldots, n$ is satisfied, and such that $depth\ I(d_1)\ge 1$, $depth\ I(d_i)\ge 2$ for all $i\ge 2$.

Let us fix a format $(f_0,\ldots ,f_n)$. We consider the pairs $(S, \GG)$ where $S$ is a Noetherian ring and $\GG$ is a free resolution of format $(f_0,\ldots ,f_n)$ over $S$.
We say that a pair $(R_{gen}, \FF_{gen})$ is a generic free resolution of format $(f_0,\ldots ,f_n)$ (or that $R_{gen}$ is a generic ring for this format) if for every pair $(S, \GG)$
there exists a (note necessarily unique) homomorphism of rings $\phi: R_{gen}\rightarrow S$ such that $\GG=\FF_{gen}\otimes_{R_{gen}}S$. The main problem considered in this paper is the existence and structure of generic rings for  formats of length three.

Throughout the paper we will deal with the notion of depth of  ideals over the rings that are not necessarily Noetherian.
This is based on the results of Northcott from his book \cite{N76}. 

 He defined the true grade
$$Grade_R (I, M)=sup_{n\ge 0} grade_{R[x_1,\ldots ,x_n]}(I\otimes_R R[x_1,\ldots ,x_n], M\otimes_R R[x_1,\ldots ,x_n]$$
and proved that  Theorems 2.1 and 2.2 hold with this definition of the depth over any ring $R$.
Thus we will use the theory over arbitrary rings and in case of (possibly) non-Noetherian ring, depth will mean the true grade in the above sense.

\bigskip\bigskip

\bigskip\bigskip

\head \S 3. Generalities on Kac-Moody Lie algebras.\endhead

\bigskip

\leftline{\bf 1. Kac-Moody Lie algebras}

\bigskip

In this section we recall the basic definitions of Kac-Moody Lie algebras. The main reference is \cite{K92}.

Let $A$ be {\it a generalized Cartan matrix}, i.e. an $n\times n$ integer matrix
$$A=(a_{i,j})_{1\le i,j\le n}$$
such that $a_{i,i}=2$ for $i=1,\ldots, n$, $a_{i,j}\le 0$ for $i\ne j$ and such that $a_{i,j}=0$ implies $a_{j,i}=0$.

We will actually assume that $A$ is {\it symmetrizable}, i.e. there exist a diagonal matrix $D$ with diagonal entries $\epsilon_1 ,\ldots ,\epsilon_n$ and a symmetric matrix $B$ such that
$$A=DB.$$

Let $l=rank (A)$. Consider the complex vector space ${\goth h}$ of dimension $2n-l$. 
We take $\Pi =\lbrace \alpha_1,\ldots ,\alpha_n\rbrace$ in ${\goth h}^*$ to be the coordinate functions.
   This is the basis of {\it simple roots}.
   We take $\Pi^\vee =\lbrace \alpha_1^\vee,\ldots ,\alpha_n^\vee\rbrace$  in ${\goth h}$ such that
   $$\langle \alpha_i^\vee, \alpha_j\rangle = a_{i,j}.$$
   Thus we can think of $\alpha_i^\vee$ as the $i$-th column of $A$. The set $\Pi^\vee$ is a basis of {\it simple coroots}.
   
   The {\it Kac-Moody Lie algebra} ${\goth g}(A)$ is the Lie algebra generated by elements $e_i, f_i$, $1\le i\le n$, with the following defining relations:
   
   $$[h, h']=0 \ for\ h,h'\in {\goth h},$$
   $$[h, e_j]=\langle h,\alpha_j\rangle e_j, [h,f_j]=-\langle h,\alpha_j\rangle f_j,$$
   $$[e_i, f_j]=\delta_{i,j}\alpha_i^\vee,$$
   $$ad (e_i)^{1-a_{i,j}}(e_j))= ad(f_i)^{1-a_{i,j}}(f_j)=0\ for\ i\ne j.$$
   
   Let $Q=\oplus_{i=1}^n \ZZ\alpha_i$,  $Q_+=\oplus_{i=1}^n \ZZ_+\alpha_i$, and $Q_-=-Q_+$.
   We define a partial ordering $\ge$ on ${\goth h}^*$ by $\lambda\ge \mu$ if and only if $\lambda-\mu\in Q_+$.
   The Kac-Moody Lie algebra ${\goth g}={\goth g}(A)$ has {\it the root space decomposition} ${\goth g}=\oplus_{\alpha\in Q}{\goth g}_\alpha$, where
   ${\goth g}_\alpha=\lbrace x\in {\goth g}\ |\ [h,x]=\alpha(h)x\rbrace$. An element $\alpha\in Q$ is called {\it a root} if $\alpha\ne 0$ and ${\goth g}_\alpha\ne 0$.
   The number $mult(\alpha)=dim\ {\goth g}_\alpha$ is called {\it the multiplicity} of the root $\alpha$. A root $\alpha>0$ (resp. $\alpha<0$) is called {\it positive} (resp. {\it negative}).
   One can easily show that every root is either positive or negative.
   We denote by $\Delta ,\Delta_+, \Delta_-$ the sets of roots, positive and negative roots respectively.
   
   For $\alpha =\sum_{i=1}^n k_i\alpha_i\in Q$ the number $ht(\alpha):=\sum_{i=1}^n k_i$ is called {\it the height} of $\alpha$.
   We define {\it the principal gradation} on ${\goth g}=\oplus_{j\in\ZZ}{\goth g}_j$ by setting ${\goth g}_j=\oplus_{ht(\alpha )=j}{\goth g}_\alpha$. Note that ${\goth g}_0={\goth h}$, ${\goth g}_{-1}=\oplus\CC f_j$,
   ${\goth g}_1=\oplus\CC e_i$. Let ${\goth g}_\pm =\oplus_{j\ge 1} {\goth g}_{\pm j}$. Then we have {\it the principal triangular decomposition}
   $${\goth g}={\goth g}_+\oplus {\goth g}_0\oplus {\goth g}_-.$$
   
   {\it The Weyl group} of $A$ is a subgroup of $Aut({\goth h}^*)$ generated by {\it the simple reflections}
   $$r_i (\lambda )=\lambda -\langle \lambda,\alpha_i^\vee\rangle\alpha_i.$$
   for $\lambda\in {\goth h}^*$.

    We choose the weight $\rho\in{\goth h}^*$ by requiring
   $$\langle\rho, \alpha_i^\vee\rangle =1,\  for \ i=1,\ldots ,n.$$
   
   A $\goth g$-module $V$ is ${\goth h}$-{\it diagonalizable} if $V=\oplus_{\lambda\in{\goth h}^*} V_\lambda$ where $V_\lambda=\lbrace v\in V\ |\ h^.v=\lambda(h)v\ \forall h\in{\goth h}\rbrace$ is the $\lambda$ {\it weight space}. If $V_\lambda\ne 0$ then $\lambda$ is a weight of $V$. The number $mult_\lambda(V):=dim\ V_\lambda$ is called {\it the multiplicity} of $\lambda$ in $V$. When all the weight spaces are finite dimensional we define {\it the character} of $V$ to be
   $$ch\ V=\sum_{\lambda\in{\goth h}^*} (dim V_\lambda) e^{\lambda},$$
   where $e^\lambda$ are the basis elements of the group algebra $\CC[{\goth h}^*]$ with the binary operation $e^\lambda e^\mu=e^{\lambda+\mu}$.
   
   Let $P(V)$ be the set of weights in $V$ and let $D(\lambda )=\lbrace \mu\in{\goth h}^*\ |\ \mu\le\lambda\rbrace $. We define the category $\OOO$ as follows: its objects are $\goth h$-diagonalizable $\goth g$-modules with finite dimensional weight spaces such that there exist finitely many elements $\mu_1,\ldots ,\mu_s$ with $P(V)\subset\cup_{i=1}^s D(\mu_i)$, and the morphisms are $\goth g$-module homomorphisms. An $\goth h$-diagonalizable module $V$ is said to be integrable if all the $e_i$, $f_i$ ($i=1,\ldots ,n$) are locally nilpotent on $V$.
   All the integrable modules in the category $\OOO$ are completely reducible (\cite{K92}, Corollary 10.7). 
   
   A ${\goth g}$-module $V$ is called {\it a highest weight module} with highest weight $\lambda$ if there is a nonzero vector $v\in V$ such that (i) ${\goth g}_+^.v=0$, (ii) $h^.v=\lambda(h)v$ for all $h\in{\goth h}$, (iii) $U({\goth g})^.v=V$. The vector $v$ is called {\it a highest weight vector}. Let ${\goth b}_+={\goth h}+{\goth g}_+$ be the Borel subalgebra of $\goth g$ and $\CC_\lambda$ the one dimensional ${\goth b}$-module defined by ${\goth g}_+ ^.1=0$, $h^.1=\lambda(h)1$ for $h\in{\goth h}$. The induced module $M(\lambda )= U({\goth g})\otimes_{U({\goth b})}\CC_\lambda$ is called {\it the Verma module} with highest weight $\lambda$. Every highest weight $\goth g$-module with highest weight $\lambda$ is a quotient of $M(\lambda )$. The Verma module contains a unique maximal proper submodule $J(\lambda )$.
   Hence the quotient $V(\lambda ):=M(\lambda )/J(\lambda )$ is irreducible, and we have a bijection between ${\goth h}^*$ and the set of irreducible modules in the category $\OOO$ given by $\lambda\mapsto V(\lambda )$.
   
   If $\lambda$ is {\it dominant integral} i.e. $\lambda (\alpha_i^\vee)\in\ZZ_+$ for all $i=1,\ldots ,n$, then $V(\lambda)$ is integrable and we have the Weyl-Kac character formula (\cite{K92}, Theorem 10.4)
   $$ch\ V(\lambda)= {{\sum_{w\in W} (-1)^{l(w)} e^{w(\lambda+\rho)-\rho}}\over{\prod_{\alpha\in\Delta_+}(1-e^{-\alpha})^{dim\ {\goth g}_\alpha}}}.$$
   Here $\rho$ is given by $\rho(\alpha^\vee_i )=1$ for $i=1,\ldots ,n$. When $\lambda=0$ we obtain the {\it denominator identity}
   $$\sum_{w\in W} (-1)^{l(w)} e^{w(\lambda+\rho)-\rho}=\prod_{\alpha\in\Delta_+}(1-e^{-\alpha})^{dim\ {\goth g}_\alpha}$$
   
   \vskip 2cm
   
  \leftline{\bf 2. Kostant formula.}
  
  \vskip 1cm
   
   Let us choose a subset $S\subset\Pi$. This defines the grading on the Kac-Moody Lie algebra
   $${\goth g}(A)=\oplus_{m\in \ZZ} {\goth g}(A)^{(S)}_m .$$
   For $m\ne 0$ the component ${\goth g}(A)^{(S)}_m$ is the span of root spaces ${\goth g}(A)_\alpha$ where $\alpha$ is a root which written in a basis of simple roots has $m$ as a  sum of coefficients of $\alpha_i$ with $\alpha_i\notin S$. Such $m$ is denoted $ht^S (\alpha)$.
   For $m=0$ ${\goth g}(A)^{(S)} _0$ also includes the Cartan subalgebra $\goth h$.
   We denote
   $${\goth g}(A)^{(S)}_+ =\oplus_{m> 0}  {\goth g}(A)^{(S)}_m , {\goth g}(A)^{S)}_- =\oplus_{m< 0}  {\goth g}(A)^{(S)}_m,$$
   so we have
   $${\goth g}(A)= {\goth g}(A)^{(S)}_+\oplus {\goth g}(A)^{(S)}_0\oplus {\goth g}(A)^{(S)}_- .$$

   We also define the subalgebra ${\goth g}_{S}$, Weyl group $W_{S}$, $\Delta_{S}$, $\Delta^\vee_{S}$ as the objects defined for the Cartan matrix
   $A' = (a_{k,l})_{k,l\in S}$. So $W_{S}$ is generated by reflections $r_k$, $k\in S$. 
   
   Define $\Delta(S)^{\pm}=\Delta^{\pm}\setminus\Delta_{S}^{\pm}$ and similarly for $\Delta(S)$.

   We also define the subset 
   $$W(S)=\lbrace w\in W\ |\ \Phi_w\subset \Delta_+ (S)\rbrace.$$
   where $\Phi_w =\lbrace \alpha\in\Delta^+ \ |\  w^{-1}(\alpha )<0\rbrace$.
   
   Let $\CC$ be a trivial ${\goth g}$-module. The homology modules $H_k ({\goth g}^{(S)}_-, \CC)$ are obtained from the complex of ${\goth g}_0^{(S)}$-modules
   $$\ldots \bigwedge^k ({\goth g}^{(S)}_- ){\buildrel{d_k}\over\rightarrow}\bigwedge^{k-1}({\goth g}_-^{(S)})\rightarrow\ldots\rightarrow \bigwedge^1 ({\goth g}^{(S)}_- ){\buildrel{d_1}\over\rightarrow}\bigwedge^{0}({\goth g}_-^{(S)}){\buildrel{d_0}\over\rightarrow}\CC\rightarrow 0,$$
   where the differential $d_k:  \bigwedge^k ({\goth g}^{(S)}_- )\rightarrow\bigwedge^{k-1}({\goth g}_-^{(S)})$ is defined as follows
   $$d_k (x_1\wedge\ldots\wedge x_k)= \sum_{s<t} (-1)^{s+t}[x_s ,x_t]\wedge x_1\ldots\wedge{\hat x}_s\wedge\ldots\wedge{\hat x}_t\wedge\ldots\wedge x_k$$
   for $k\ge 2$, $x_i\in {\goth g}^{(S)}_-$, and $d_1=d_0=0$.
   
   For simplicity we write $H_k ({\goth g}^{(S)}_-)$ instead of $H_k ({\goth g}^{(S)}_-, \CC)$. Each of the terms $ \bigwedge^k ({\goth g}^{(S)}_- )$ has a $\ZZ$-grading induced by that on ${\goth g}^{(S)}_-$.
   For $j\ge 0$ we define $\bigwedge^k ({\goth g}^{(S)}_-)_{-j}$ to be the subspace of $ \bigwedge^k ({\goth g}^{(S)}_- )$ spanned by the vectors of the form $x_1\wedge\ldots\wedge x_k$
   such that $deg(x_1)+\ldots +deg (x_k)=-j$. The homology module $H_k ({\goth g}^{(S)}_-)$ also has the induced $\ZZ$-grading. Note that $\bigwedge^k ({\goth g}^{(S)}_-)_{-j}=H_k ({\goth g}^{(S)}_-)_{-j}=0$
   for $k>j$. The ${\goth g}_0^{(S)}$-module structure of the homology modules $H_k ({\goth g}^{(S)}_-)$ is determined by the following formula known as {\it Kostant's formula}.
   
   \proclaim{Theorem 3.1.(\cite{GL76}, \cite{Li92})}
   $$H_k ({\goth g}^{(S)}_-)=\oplus_{w\in W(S), l(w)=k} V_S(w\rho -\rho),$$
   where $V_S(\lambda)$ denotes the integrable highest weight ${\goth g}_0^{(S)}$-module with highest weight $\lambda$.
   \endproclaim

   \vskip 2cm

   \vskip 2cm
   
   \leftline{\bf 3. Parabolic version of a BGG resolution.} 
   
   \vskip 1cm
   
   \proclaim{Definition 3.2} Let $S\subset \lbrace 1,\ldots ,n\rbrace$.
   Let ${\goth p}^{(S)}=\oplus_{j\ge 0} ({\goth g}^{(S)})_{-j}$ be the parabolic subalgebra.
   Define a generalized Verma module
   $$M(\lambda )^{(S)}= U({\goth g})\otimes_{U({\goth p}^{(S)})}L_S(\lambda)$$
   where $V_S(\lambda )$ is considered as a ${\goth g}^{(S)}_0\oplus {\goth g}^{(S)}_+$-module where ${\goth g}^{(S)}_+$ acts trivially.
   \endproclaim
   
   \proclaim{Theorem 3.3.(\cite{Ku02}, section 9.2)}
   There exists an exact complex of ${\goth p}^{(S)}$-\break modules
   $$\ldots \rightarrow F^p_{(S)}\rightarrow\ldots \rightarrow F^1_{(S)}\rightarrow F^0_{(S)}\rightarrow V^{}(\lambda )\rightarrow 0$$
   where
   $$F^p_{(S)}=\oplus_{w\in W'_{(S)}, l(w)=p} M(w^{-1}{}^.\lambda)^{(S)}.$$
   Here $V(\lambda)$ is the irreducible ${\goth g}$-module of highest weight $\lambda$, $w^.\lambda :=w(\lambda+\rho )-\rho$, and 
   $$W'_{(S)}=\lbrace w\in W\ |\ l(wv)\ge l(w)\ \forall v\in W_{(S)}\rbrace.$$
   where $W_{(S)}$ denotes the subgroup of $W$ generated by $r_i$ ($i\in S$).
   Thus $W'_{(S)}$ is the set of elements of minimal length in the cosets of $W_{(S)}$ (there is one in each coset).
    \endproclaim
   
   \vskip 1cm

\head  \S 4. The Kac-Moody Lie algebra of type $T_{p,q,r}$ \endhead

\bigskip

We will be interested in the special diagrams $T_{p,q,r}$ defined as follows
$$\matrix x_{p-1}&-&x_{p-2}&\ldots&x_{1}&-&u&-&y_{1}&\ldots&y_{q-2}&-&y_{q-1}\\
&&&&&&|&&&&&&\\
&&&&&&z_{1}&&&&&&\\
&&&&&&|&&&&&&\\
&&&&&&\ldots&&&&&&\\
&&&&&&z_{r-2}&&&&&&\\
&&&&&&|&&&&&&\\
&&&&&&z_{r-1}&&&&&&
\endmatrix$$

We recall some basic notions about Kac-Moody Lie algebra ${\goth g}(T_{p.q.r})$ associated to this diagram.
The generalized Cartan matrix $A(T_{p,q,r})$ has rows and columns indexed by the set
$\lbrace 0,1,\ldots, p-1, 1',\ldots ,(q-1)', 1'',\ldots ,(r-1)''\rbrace$ corresponding to the vertices $u, x_1 ,\ldots ,x_{p-1}, y_1,\ldots ,y_{q-1}, z_1,\ldots ,z_{r-1}$ respectively.
Sometimes we denote vertices by natural numbers from $[1,p+q+r-2]$, in the order listed above.

The entries of $A$ are given by
$$A(T_{p,q,r})_{i,j} =\cases 2&\ if\ i=j;\\
                                                  -1&\  if\ the\ nodes\ i\ and\  j\  are\ incident\ in\  T_{p,q,r};\\
                                                  0\ & otherwise.
                                                  \endcases$$
                                                  
   We set $n:=p+q+r-2$ so $A(T_{p,q,r})$ is an $n\times n$ matrix.
   The following is an easy consequence of results in  \cite{K92}
   
   \proclaim{Proposition 4.1} 
   \item{a)} If $T_{p,q,r}$ is a Dynkin diagram, then  the matrix $A(T_{p,q,r})$ has rank  $n$. The quadratic form
   corresponding to $A(T_{p,q,r})$ is positive definite,
    \item{b)} If $T_{p,q,r}$ is an affine Dynkin diagram, then  the matrix $A(T_{p,q,r})$ has rank  $n-1$. The quadratic form
   corresponding to $A(T_{p,q,r})$ is semi-positive definite,
    \item{a)} In all other cases  the matrix $A(T_{p,q,r})$ has rank  $n$. The quadratic form
   corresponding to $A(T_{p,q,r})$ has signature $(n-1,1)$,
   \endproclaim                            
   
   \demo {Proof} The first two statements are special cases of Theorem 4.3, the last is exercise 4.6 from \cite{K92}.
   \qqed
   \enddemo
   
   Let us describe the roots, coroots and the Weyl group. We take the vector space $\goth h$ of dimension $n$ if $T_{p,q,r}$ is not affine and $n+1$ if it is.

   Assume first that $T_{p,q,r}$ is not affine. We take $\Pi =\lbrace \alpha_1,\ldots ,\alpha_n\rbrace$ in ${\goth h}^*$ to be the coordinate functions.
   This is the basis of simple roots.
   We take $\Pi^\vee =\lbrace \alpha_1^\vee,\ldots ,\alpha_n^\vee\rbrace$  in ${\goth h}$ such that
   $$\langle \alpha_i^\vee, \alpha_j\rangle = a_{i,j}.$$
   Thus we can think of $\alpha_i^\vee$ as the $i$-th column of $A(T_{p,q,r})$. The set $\Pi^\vee$ is a basis of simple coroots.
   
     As defined in the previous section, the  Weyl group of $A$ is a subgroup of $Aut({\goth h}^*)$ generated by the simple reflections.
   $$r_i (\lambda )=\lambda -\langle \lambda,\alpha_i^\vee\rangle\alpha_i.\eqno{(1)}$$
   for $\lambda\in {\goth h}^*$. For the graph $T_{p,q,r}$ (or any tree) there is also a combinatorial formula
   
   $$r_i (\lambda_1,\ldots,\lambda_n )= (\lambda_1,\ldots ,\lambda_n )-2\lambda_i\alpha_i  +\sum_{x-i}\lambda_i\alpha_x .\eqno{(2)}$$
   This means that to calculate the value of the reflection $r_i$ on $\lambda$ (thought of as a graph $T_{p,q,r}$ with labeled vertices), we reverse the sign of a label at $i$,
   and add this label to the labels of all neighbors of $i$.
   
   For the affine $T_{p, q, r}$ the above formulas are still true (one has to remember that $\alpha_1^\vee ,\ldots ,\alpha_n^\vee$ are still independent because they also have a component
   on the $(n+1)$'st coordinate $\alpha_{n+1}$ which is not a part of a basis of simple roots).

   We now specialize to $A=A(T_{p,q,r})$ and to $S=[1,n]\setminus \lbrace p+q\rbrace$. This means the distinguished root is the root corresponding to the vertex $z_1$.
   We will write ${\goth g}:={\goth g}(T_{p,q,r})$ and ${\goth g}_i:={\goth g}^{(S)}_i$  in this case. We have the following proposition.
   
   \proclaim{Proposition 4.2} We have
\item{a)} ${\goth g}_0= {\goth sl}_{r-1}\times{\goth sl}_{p+q}\times\CC$,
\item{b)} ${\goth g}_1= \CC^{r-1}\otimes \bigwedge^{p}\CC^{p+q}$,
\item{c)} 
$${\goth g}_2 = \bigwedge^2\CC^{r-1}\otimes Ker(S_2 (\bigwedge^p \CC^{p+q})\rightarrow S_{2^p}\CC^{p+q})\oplus $$
$$\oplus S_2\CC^{r-1}\otimes Ker(\bigwedge^2 (\bigwedge^p\CC^{p+q})\rightarrow S_{2^{p-1},1^2}\CC^{p+q}).$$
\item{d)} The higher components $\LL_m$ can be defined as cokernels of the graded components of the Koszul complex
$$(\bigwedge^3 \LL)_m\rightarrow (\bigwedge^2\LL)_m\rightarrow \LL_m\rightarrow 0.$$
   \endproclaim

\demo{Proof} We use the generalized Kostant formula to identify ${\goth g}(T_{p,q,r})^{(S)}_{>0}$ for $S=[1,n]\setminus\lbrace p+q\rbrace$.
We denote by $s_i$ the simple reflection corresponding to the vertex $i$, where vertices are labeled by $0$, $1,\ldots ,p-1$, $1',\ldots ,(q-1)'$, $1'',\ldots ,(r-1)''$ as in section 3.
The only elements of length two in the subgroup $W(S)$ are  the elements $s_{1''}s_0$ and $s_{1''}s_{2''}$. We identify a weight with a labeled Dynkin diagram, each vertex labeled by a coefficient of the coresponding fundamental weight.
Calculating the corresponding values
of $w\rho -\rho$ and using the formula $(2)$  we get
$$ s_{1''}s_0\rho -\rho =\matrix 0&-&0&\ldots&1&-&0&-&1&\ldots&0&-&0\\
&&&&&&|&&&&&&\\
&&&&&&-3&&&&&&\\
&&&&&&|&&&&&&\\
&&&&&&2&&&&&&\\
&&&&&&|&&&&&&\\
&&&&&&\ldots&&&&&&\\
&&&&&&0&&&&&&\\
&&&&&&|&&&&&&\\
&&&&&&0&&&&&&
\endmatrix$$
with values zero at all other vertices. This (discarding labeling of the vertex $1''$) corresponds to the representation $S_2\CC^{r-1}\otimes  S_{2^{p-1},1^2}\CC^{p+q}$.
$$ s_{1''}s_{2''}\rho-\rho=\matrix 0&-&0&\ldots&0&-&2&-&0&\ldots&0&-&0\\
&&&&&&|&&&&&&\\
&&&&&&-3&&&&&&\\
&&&&&&|&&&&&&\\
&&&&&&0&&&&&&\\
&&&&&&|&&&&&&\\
&&&&&&1&&&&&&\\
&&&&&&|&&&&&&\\
&&&&&&\ldots&&&&&&\\
&&&&&&0&&&&&&\\
&&&&&&|&&&&&&\\
&&&&&&0&&&&&&
\endmatrix$$
with values zero at all other vertices. This (discarding labeling of the vertex $1''$) corresponds to the representation $\bigwedge^2\CC^{r-1}\otimes  S_{2^p}\CC^{p+q}$.
Since these are the only weights in $H_2 ({\goth g}_+^{(S)})$,  description d) follows.
\qqed
\enddemo

\bigskip\bigskip

Taking $E=\CC^{r-1}$, $F=\CC^{p+q}$ we denote $\LL(p, E, F)$ the positive part
$$\LL(p,E, F)=\oplus_{i>0}{\goth g}_i .$$

\proclaim {Proposition 4.3} The algebra $\LL (p, E, F)$ is finite dimensional if and only if one of the following cases occurs.
\item{a)} $p=q=2$, $r\ge 2$ arbitrary,
\item{b)}$q=r=2$, $p\ge 3$ arbitrary,
\item{c)}$q=2$, $r=3$, $p=3,4,5$,
\item{d)} $q=3$, $r=2$, $p=3,4,5$,
\item{e)} $q=2$, $p=3$, $r=4,5$.
\endproclaim

\demo{Proof} 
Indeed, the listed cases are exactly the cases when $T_{p,q,r}$ is a Dynkin diagram. To be more precise, we have $T_{p,q,r}=D_{r+2}$ in case $a)$,
$T_{p,q,r}=D_{p+2}$ in case $b)$, $T_{p,q,r}= E_{p+3}$ in case $c)$, $T_{p,q,r}=E_{p+3}$ in case $d)$ and $T_{p,q,r}=E_{r+3}$ in case $e)$.
In the listed cases the algebra $\LL (p,E, F)$ is obviously finite dimensional.
In the other cases, the positive part of the Kac-Moody Lie algebra is infinite dimensional.
\qqed
\enddemo

Let us also calculate the beginning part of the parabolic BGG resolutions for the case under consideration, i.e. ${\goth g}:={\goth g}(T_{p,q,r})$,
$S=\lbrace [1,n]\setminus \lbrace p+q\rbrace\rbrace$.

\proclaim{Proposition 4.4} Let ${\goth g}:={\goth g}(T_{p,q,r})$,
$S=\lbrace [1,n]\setminus \lbrace p+q\rbrace\rbrace$.
Let us consider the highest weight $\lambda =(\lambda_1,\ldots ,\lambda_n)$. 
The three initial terms of the parabolic BGG complex described in Theorem 3.3 are.
$$\matrix M(r_{p+q}r_{1}{}^.\lambda )^{(S)}\oplus M(r_{p+q}r_{p+q+1}{}^.\lambda )^{(S)}\\
\downarrow\\
M(r_{p+q}{}^.\lambda)^{(S)}\\
\downarrow\\
M(\lambda )^{(S)}
\endmatrix$$
where $M(\mu )^{(S)}$ denotes the parabolic Verma module.

To make things explicit we identify the weights
$$r_{p+q}{}^.\lambda = (\lambda_1+\lambda_{p+q}+1,\ldots ,-\lambda_{p+q}-2, \lambda_{p+q}+\lambda_{p+q+1}+1,\ldots ),$$
where listed components are at vertices $1, p+q, p+q+1$ and not listed ones are the same as in $\lambda$.
$$r_{p+q}r_{1}{}^.\lambda = $$
$$(\lambda_{p+q}-1, \lambda_1+\lambda_{2}+1,\ldots, \lambda_1+\lambda_{p+1}+1, \ldots ,-\lambda_1-\lambda_{p+q}-3, \lambda_1+\lambda_{p+q}+\lambda_{p+q+1}+1,\ldots )$$
where listed components are at vertices $1, 2, p, p+q, p+q+1$ and not listed ones are the same as in $\lambda$.
$$r_{p+q}r_{p+q+1}{}^.\lambda = $$
$$(\lambda_1+\lambda_{p+q}+\lambda_{p+q+1}+1,\ldots ,-\lambda_{p+q}-\lambda_{p+q+1}-5, \lambda_{p+q}-1, \lambda_{p+q+1}+\lambda_{p+q+2}+1,\ldots )$$
where listed components are at vertices $1, p+q, p+q+1, p+q+2$ and not listed ones are the same as in $\lambda$.
\endproclaim

\bigskip\bigskip

\head \S 5. Geometric approach to acyclicity.\endhead

We prove the geometric result on acyclicity of free complexes. It is based on homological algebra and it is essential for our approach.

\proclaim {Theorem 5.1}   Let $X=Spec R$, and let $j:U\rightarrow X$ be an open immersion. Let
   $$\GG:0\rightarrow G_n\rightarrow G_{n-1}\rightarrow \ldots\rightarrow G_1\rightarrow G_0$$
   be a complex of free $R$-modules (treated as a complex of sheaves over $X$)  such that $\GG\ |_U$ is acyclic. Then $H_n (\GG\otimes j_*\OOO_U)=0$, 
    $H_{n-1} (\GG\otimes j_*\OOO_U)=0$, and the complex $\GG\otimes j_*\OOO_U$ is acyclic if and only if
    $\RRR^i j_*\OOO_U=0$ for $i=1,\ldots , n-2$. \endproclaim

\demo{Proof} 
Before we start, let us decompose the complex $\GG\ |_U$ to short exact sequences.
Denoting $B_i = Im(G_{i+1}\ |_U\rightarrow G_i\ |_U)$ we have
exact sequences 
$$0\rightarrow B_{i}\rightarrow G_i\ |_U\rightarrow B_{i-1}\rightarrow 0$$
for $i=2,\ldots, n$ (with $B_{n}=0$). and
$$0\rightarrow B_1\rightarrow G_1\ |_U\rightarrow G_0\ |_U.$$
This induces long exact sequences
$$0\rightarrow j_* B_{i}\rightarrow G_i\otimes j_*\OOO_U\rightarrow j_*B_{i-1}\rightarrow R^1j_* B_{i}\ldots$$
as well as an exact sequence
$$0\rightarrow j_*B_1\rightarrow G_1\otimes j_*\OOO_U\rightarrow G_0\otimes j_*\OOO_U .$$

Next we show that vanishing of higher direct images implies acyclicity.
Indeed, our vanishing implies that
$R^i j_* B_{n-s}=0$ for $1\le i\le n-s-1$. So the above exact sequences have last term zero and we get that $\GG\otimes j_*\OOO_U$ is acyclic.

To prove the reverse implication let us proceed by induction on $n$.
For $n=2$ there is nothing to prove. For $n=3$ we see from the exact sequences that $H_3(\GG)=H_2(\GG)=0$ and
$H_1(\GG)=Ker (R^1j_*G_3\rightarrow R^1j_*G_2)$. We need

\proclaim{Lemma 5.2} Let $M$ be an $R$-module. Let $\phi:F\rightarrow G$ be  a map of free $R$-modules of finite rank. Denote by $I(\phi)$ the ideal of maximal minors of $\phi$. Then
$\phi\otimes M$ is a monomorphism if and only if $depth_R(I(\phi), M)\ge 1$.
\endproclaim

\demo{Proof} This is a special case of Theorem 2, Appendix B from \cite{N76}.
\enddemo

In our case $R^1j_*{\OOO_U}$ is supported on $X\setminus U$ so $Ker (R^1j_*G_3\rightarrow R^1j_*G_2)=0$ implies $R^1j_*\OOO_U=0$, completing the case $n=3$.
Assume the result is proved for $n-1$ and the complex $j_*\GG=\GG\otimes j_*\OOO_U$ is acyclic. By induction (applied to $\GG$ truncated at $G_1$) we have
$$R^1j_*\OOO_U=R^2j_*\OOO_U=\ldots =R^{n-3}J_*\OOO_U=0.$$
Now our long exact sequences imply that $R^{n-3}j_* B_{n-2}=R^{n-4}j_* B_{n-3}=\ldots =R^1j_* B_2.$
We also have the exact sequence
$$0\rightarrow R^{n-3}j_* B_{n-2}\rightarrow R^{n-2}j_* G_n\rightarrow R^{n-2}j_* G_{n-1}.$$
Also, from the exact sequences we can deduce that 
$$H_1(\GG\otimes j_*\OOO_U)= Ker(R^1j_* B_2\rightarrow R^1j_* G_2)=R^1j_* B_2$$
This means that if $H_1(\GG\otimes j_*\OOO_U)=0$  then $R^1j_* B_2=0$, so the map
$ R^{n-2}j_* G_n\rightarrow R^{n-2}j_* G_{n-1}.$ is a monomorphism, which implies by  Lemma 5.2 that
$R^{n-2}j_*\OOO_U=0$. 
\qqed
\enddemo

\head \S6. The rings $R_a$ generated by Buchsbaum-Eisenbud multipliers. \endhead

\bigskip

In this section we recall properties of the rings $R_a$ which are obtained from coordinate rings of the 
varieties of generic complexes by adding the 
Buchsbaum-Eisenbud multipliers and factoring the relations satisfied by them. These rings are the starting point 
of our construction. Their properties (rational singularities and sphericality) are essential for the whole approach.
Most of the results of this section were proved in \cite{PW90}, section 1. The additional results are easy consequences.

In what follows we use heavily representation theory of $GL_n$. For a $GL_n$ dominant weight $\lambda = (\lambda_1,\ldots ,\lambda_n)\in \ZZ^n$ we denote
$S_\lambda F=S_{(\lambda_1, \lambda_2,\ldots ,\lambda_n)} F$ {\it the Schur functor} on the space $F=K^n$.

Let us fix the format $(f_0,\ldots ,f_n)$. We work over a fixed field $K$. 
In this section we assume that $K$ has characteristic zero.
Consider the variety $X_c$ of complexes
$$0\rightarrow F_n{\buildrel {d_n}\over\rightarrow}F_{n-1}\rightarrow\ldots\rightarrow F_1{\buildrel{d_1}\over\rightarrow} F_0$$
of vector spaces over $K$, with $rank\ F_i=f_i$ and $rank\ d_i\le r_i$. We fix bases $\lbrace e^i_{j(i)}\rbrace_{1\le j(i)\le f_i}$ of $F_i$ for each $i=1,\ldots ,n$.

The coordinate ring $R_c$ of $X_c$ can be obtained as follows. We add to $K$ the variables $X^{(i)}_{j(i),k(i)}$ ($1\le i\le n, 1\le j(i)\le f_{i-1}, 1\le k(i)\le f_i$)
which are the entries of the generic maps $d_i$ in bases $\lbrace e^i_{j(i)}\rbrace_{1\le j(i)\le f_i}$. The corresponding $f_{i-1}\times f_i$ matrix of variables over $S_c$ is denoted $X^{(i)}$.
We denote the resulting polynomial ring $S_c$.
We define the ideal $J_c$ in the polynomial ring $S_c$ as follows. $J_c$ is generated by the entries of matrices $X^{(i-1)}\circ X^{(i)}$ ($2\le i\le n$) and by the $(r_i+1)\times (r_i+1)$ minors of $X^{(i)}$ ($1\le i\le n-1$).
Finally we define $R_c:=S_c/J_c$. The ring $R_c$ is the coordinate ring of $X_c$ which is {\it the variety of generic complexes}.

Recall  that the variety $X_c$ has a natural desingularization $Z_c$. For $i=1,\ldots ,n-1$ denote by $Grass(r_{i+1}, F_i)$ the Grassmannian of subspaces of rank $r_{i+1}$ of $F_i$.
Let 
$$0\rightarrow {\RRR}_i\rightarrow F_i\times Grass (r_{i+1}, F_i)\rightarrow {\QQQ}_i\rightarrow 0$$
be the  tautological sequence on $Grass(r_{i+1}, F_i)$.

$$Z_c=\lbrace ((d_1,\ldots ,d_n),(R_{1},\ldots ,R_{n-1}))\in X_c\times\prod_{i=1}^{n-1} Grass(r_{i+1} ,F_{i})\ |\ Im(d_{i+1})\subset R_i\ \rbrace .$$
Denote $p:Z_c\rightarrow \prod_{i-1}^{n-1}$, $q:Z_c\rightarrow X_c$ the natural projections.
We have $p_* {\OOO}_{Z_c}=\otimes_{i=1}^n Sym({\QQQ_i\otimes \RRR}^*_{i-1})$, where ${\QQQ}_n=F_n$.

\proclaim{Theorem 6.1} (\cite{DCS81}, \cite{PW90}) The variety $X_c$ carries the natural action of the group $GL:=\prod_{i=0}^n GL(F_i)$.
It is a spherical variety and it has rational singularities.
The coordinate ring $R_c$ has a multiplicity free decomposition to the irreducible representations of $GL$ given by the formula
$$R_c=\oplus_{\alpha^{(1)},\ldots ,\alpha^{(n)}} \otimes_{i=0}^n S_{(\alpha^{(i)}_1,\ldots ,\alpha^{(i)}_{r_i}, -\alpha^{(i+1)}_{r_{i+1}},\ldots ,-\alpha^{(i+1)}_1)}F_i$$
where we sum over all $n$-tuples $(\alpha^{(1)},\ldots ,\alpha^{(n)})$ of partitions, with the $i$-th partition $\alpha^{(i)}=(\alpha^{(i)}_1,\ldots \alpha^{(i)}_{r_i})$ having at most $r_i$ parts.
Here by convention $\alpha^{(n+1)}=0$ has no parts and $\alpha^{(0)} =(0^{f_0-r_1})$ has $f_0-r_1$ parts.
\endproclaim

This result follows by standard methods from Cauchy decomposition
$$\otimes_{i=1}^n Sym({\QQQ_i\otimes \RRR}^*_{i-1}) = \oplus_{\alpha^{(1)},\ldots ,\alpha^{(n)}} \otimes_{i=1}^n S_{\alpha^{(i)}}\QQQ_i\otimes S_{\alpha^{(i)}}\RRR^*_{i-1},$$
and the fact that by Bott theorem the higher cohomology of the above sheaf vanishes and the sections decompose as given in the Theorem 6.1.

We have {\it a generic complex $\FF_c$ of format $(f_0,\ldots ,f_n)$} defined over the ring $R_c$. It is a complex
$$\FF_c:0\rightarrow F_n\otimes R_c{\buildrel{d_n}\over\rightarrow}F_{n-1}\otimes R_c\rightarrow\ldots\rightarrow F_1\otimes R_c{\buildrel{d_1}\over\rightarrow}F_{0}\otimes R_c$$
with $d_i$ given (in our bases of $F_i$) by the matrix $X^{(i)}$.

In \cite{PW90}, section 1 we carried a similar procedure for the rings $R_a$. 

Consider the affine space $X=\prod_{i=1}^n Hom_K (F_i, F_{i-1})\times\prod_{i=1}^{n-1}\bigwedge^{r_i}F_{i-1}$. The coordinates in $Hom_K (F_i, F_{i-1}$ are the entries of  the map $d_i$ and the coordinates of $\bigwedge^{f_i}F_{i-1}$ are the Buchsbaum-Eisenbud multipliers $a_i$.
We also consider the  analogue of the desingularization $Z_c$.
$$Z_a\subset X\times \prod_{i=1}^{n-1} Grass(r_{i+1}, F_i)$$

\proclaim{Definition 6.2}  The variety $Z_a$,
$$Z_a\subset X\times \prod_{i=1}^{n-1} Grass(r_{i+1}, F_i)$$
is defined by the following conditions.
A point 
$$\lbrace (d_1,\ldots, d_n ;a_1,\ldots ,a_{n-1}), (R_1,\ldots ,R_{n-1})\rbrace\in X\times \prod_{i=1}^{n-1} Grass(r_{i+1}, F_i)$$
 is in $Z_a$
if and only if the following conditions are satisfied.
\item{(1)} $Im(d_i)\subset R_i\subset Ker(d_{i-1}),$
\item{(2)} $a_i\in\bigwedge^i R_{i-1}$,
\item{(3)} For the induced map $d'_i: Q_i\rightarrow R_{i-1}$ we have $d'_n=a_n$, $d'_i=a_{i+1}a_i$ for $i=1,\ldots ,n-1$.
We denote $p:Z_a\rightarrow X$, $q:Z_a\rightarrow \prod_{i=1}^{n-1} Grass(r_{i+1}, F_i)$ two projections, and we define the variety $X_a:=p(Z_a)\subset X$.
\endproclaim

The variety $Z_a$ is fiber bundle over $\prod_{i=1}^{n-1} Grass(r_{i+1}, F_i)$, and the fibre over a point $(R_1,\ldots ,R_{n-1})$  is the affine variety
given by the general $r_i\times r_i$ matrices $d'_i:Q_i\rightarrow R_{i-1}$ and elements $a_i\in\bigwedge^{r_i}R_{i-1}$ satisfying relations $(3)$.

It turns out that $Z_a$ has rational singularities so it can be used in a similar way to $Z_c$.
The following result is proved in \cite{PW90}, section 1.

\proclaim{Theorem 6.3} The variety $X_a$ carries the natural action of the group $GL:=\prod_{i=0}^n GL(F_i)$.
It is a spherical variety and it has rational singularities.
The coordinate ring $R_a$ of $X_a$ has a multiplicity free decomposition to the irreducible representations of $GL$ given by the formula
$$R_a=\oplus_{\alpha^{(1)},\ldots ,\alpha^{(n),x^{(1)},\ldots ,x^{(n)}}}$$
$$ \otimes_{i=0}^n S_{(\chi^{(i)}+\alpha^{(i)}_1,\ldots ,\chi^{(i)}+\alpha^{(i)}_{r_i-1},\chi^{(i)}, -\chi^{(i+1)}, -\chi^{(i+1)}-\alpha^{(i+1)}_{r_{i+1}-1},\ldots ,-\chi^{(i+1)}-\alpha^{(i+1)}_1)}F_i.$$
Here the notation is as follows. We sum over all $n$-tuples $(\alpha^{(1)},\ldots ,\alpha^{(n)})$ of partitions, with the $i$-th partition $\alpha^{(i)}=(\alpha^{(i)}_1,\ldots \alpha^{(i)}_{r_i-1})$ having at most $r_i-1$ parts. Here by convention $\alpha^{(n+1)}=0$ has no parts and $\alpha^{(0)} =(0^{f_0-r_1})$ has $f_0-r_1$ parts.

We also sum over $n$-tuples of natural numbers $x^{(1)},\ldots ,x^{(n)}$ (degrees with respect to the $a_i$'s).
The numbers $\chi^{(i)}$ are partial Euler characteristics and they are
$$\chi^{(i)}=\sum_{j=1}^i (-1)^{i-j} x^{(j)}.$$
Note that $\chi^{(i)}+\chi^{(i+1)}= x^{(i+1)}$, so all the weights listed in the formula are dominant.
\endproclaim

The defining relations of the ring $R_a$ are written down explicitly in \cite{PW90}, section 1. The structure of $R_a$ can be described also in a characteristic free way
(replacing $K$ by $\ZZ$ and using filtrations instead of direct sums).  This was done in \cite{PW90}, section 1 and in \cite{T01} where some errors in characteristic free part of the approach were fixed.

We denote by $\FF_a$ the complex $\FF_c\otimes_{R_c}R_a$. This complex has a weaker universality property, true even in a characteristic free version.

\proclaim{Theorem 6.4}The complex $\FF_a$ is the universal complex of format $(f_0,\ldots ,f_n)$ which is acyclic in codimension $1$.
This means that for every pair $(S, \GG)$ such that $S$ is a Noetherian ring and $\GG$ is a complex of free modules of format $(f_0,\ldots ,f_n)$ over $S$ which is acyclic of codimension $1$ (i.e. the set of points where the complex is not acyclic has a defining ideal of depth $\ge 2$), then there is a unique homomorphism
$\phi: R_a\rightarrow S$ such that $\GG=\FF_a\otimes_{R_a}S$.
\endproclaim

This has the following consequence (which goes back to Huneke and Hochster (\cite{H75}) and is even true over $\ZZ$).

\proclaim{Corollary 6.5} For $n=2$ the pair $(R_a, \FF_a)$ is a generic acyclic complex for the format $(f_0, f_1, f_2)$.
\endproclaim

In the remainder of this section we look more closely at the homology modules of the complex $\FF_a$. These modules are possible to analyze thanks to the multiplicity free
structure of the ring $R_a$.

We look at the module $F_j\otimes R_a$ and compare it to the modules $F_{j+1}\otimes R_a$ and $F_{j-1}\otimes R_a$. We  describe the cancellations that occur between them when applying the map $d_{j+1}$ and $d_j$.
Looking at the representations $F_j\otimes R_a$ and $F_{j-1}\otimes R_a$, let us assume that they have common representations coming from summands

$$F_j\otimes  \otimes_{i=0}^n S_{(\chi^{(i)}+\alpha^{(i)}_1,\ldots ,\chi^{(i)}+\alpha^{(i)}_{r_i-1},\chi^{(i)}, -\chi^{(i+1)}, -\chi^{(i+1)}-\alpha^{(i+1)}_{r_{i+1}-1},\ldots ,-\chi^{(i+1)}-\alpha^{(i+1)}_1)}F_i$$

$$F_{j-1}\otimes  _{i=0}^n S_{(\psi^{(i)}+\beta^{(i)}_1,\ldots ,\psi^{(i)}+\beta^{(i)}_{r_i-1},\psi^{(i)}, -\psi^{(i+1)}, -\psi^{(i+1)}-\beta^{(i+1)}_{r_{i+1}-1},\ldots ,-\psi^{(i+1)}-\beta^{(i+1)}_1)}F_i$$
We denote the degrees of $a_i$'s by $x^{(i)}$ in the first representation, and by $y^{(i)}$ in the second one.

We note that by Pieri formula we have two possibilities.
On the coordinate $F_j$ we can add a box to first $r_j-1$ entries, or to the $r_j$-th entry.

It is easy to see that when we add a box in $F_j\otimes R_a$ to one of the first $r_j-1$ places, we can always find a corresponding representation in $F_{j-1}\otimes R_a$, with
$\alpha^{(i)}=\beta^{(i)}$, $\chi^{(i)}=\psi^{(i)}$ for $i=1,\ldots ,n$.

Consider the crucial case when we add a box in $F_j\otimes R_a$ to the $r_j$-th place, and we add a box in $F_{j-1}\otimes R_a$ in the $(r_{j-1}+1)$'st place.
Then comparing the weights we have:
$$\chi^{(i)}=\psi^{(i)}, \alpha^{(i)}=\beta^{(i)}, \forall\ i\ne j,$$
$$\chi^{(j)}+1=\psi^{(j)}, \alpha^{(j)}_k=\beta^{(j)}_k-1,\forall\ 1\le k\le r_j-1.$$

This translates to $x^{(i)}=y^{(i)}\ \forall i\ne j,j+1,$ and $x^{(j)}=y^{(j)}-1, x^{(j+1)}=y^{(j+1)}-1$.

This means such cancellation cannot occur when $y^{(j+1)}=0$, so the corresponding representations stay in $H_{j-1}(\FF_a)$. Notice that this does not happen for $j-1=n, n-1$, as there is no $y^{(j+1)}$ in such cases.

\proclaim{Theorem 6.6} The homology groups $H_n(\FF_a)=H_{n-1}(\FF_a)=0$.
For $1\le j-1\le n-2$ we have 
$$H_{j-1}(\FF_a)= \oplus_{\beta^{(1)}, \ldots ,\beta^{(n)}, y^{(1)},\ldots ,y^{(n)}, y^{(j+1)}=0}$$
$$ \otimes  _{i=0}^{j-2} S_{(\psi^{(i)}+\beta^{(i)}_1,\ldots ,\psi^{(i)}+\beta^{(i)}_{r_i-1},\psi^{(i)}, \psi^{(i+1)}, -\psi^{(i+1)}-\beta^{(i+1)}_{r_{i+1}-1},\ldots ,-\psi^{(i+1)}-\beta^{(ji+1}_1)}F_{i}\otimes$$
$$ \otimes S_{(\psi^{(j-1)}+\beta^{(j-1)}_1,\ldots ,\psi^{(j-1)}+\beta^{(j-1)}_{r_{j-1}-1},\psi^{(j-1)}, 1-\psi^{(j)}, -\psi^{(j)}-\beta^{(j)}_{r_{j}-1},\ldots ,-\psi^{(j)}-\beta^{(j)}_1)}F_{j-1}\otimes$$
$$ \otimes  _{i=j}^n S_{(\psi^{(i)}+\beta^{(i)}_1,\ldots ,\psi^{(i)}+\beta^{(i)}_{r_i-1},\psi^{(i)}, -\psi^{(i+1)}, -\psi^{(i+1)}-\beta^{(i+1)}_{r_{i+1}-1},\ldots ,-\psi^{(i+1)}-\beta^{(i+1)}_1)}F_i.$$
In particular the homology group $H_{j-1}(\FF_a)$ is annihilated by the ideal $I(a_{j+1})$ generated by the entries of the $(j+1)$'st Buchsbaum-Eisenbud multiplier map.
\endproclaim

\demo{Proof} To prove the theorem it is enough to show that the indicated cancellations indeed occur. This is not difficult since for each $i$ the module $F_i\otimes R_a$ is multiplicity free as a $\GL(\FF)$-module. Moreover the highest weight vectors of irreducible representations are not difficult to write down as only the Pieri formula of multiplying by $F_i$ is involved. We skip the details here since the result is not used elsewhere.
\qed
\enddemo

Let us look at the generator of $H_{j-1}(\FF_a)$ for $n-2\ge j-1\ge 1$. It will be a minimal partition in $H_{j-1}(\FF_a)$.
We obtain it by setting $\beta^{(i)}=0$ for $i=1,\ldots ,n$, $\psi^{(i)}=0$ for $i>j-1$, $\psi^{(i)}=(-1)^{j-1-i}$ for $1\le i\le j-1$.
The resulting representation is (up to maximal exterior powers of $F_i$) $\bigwedge^{r_{j-1}+1}F_{j-1}$.
The existence of such cycle $q^{(j-1)}_1$ means in the generic ring one will need to add a representation $F_j^*\otimes\bigwedge^{r_{j-1}+1}F_{j-1}$, corresponding to the map
$p_1^{(j-1)}: \bigwedge^{r_{j-1}+1}F_{j-1}\otimes R\rightarrow F_j\otimes R$ covering this cycle. The maps $p_1^{(j-1)}$ are the same as the maps $b_{j-1}$ coming from the Second Structure Theorem of Buchsbaum and Eisenbud (\cite{BE74}, section 6).

Finally we note the key property of the lattice of weights of $R_a$.

\proclaim{Remark 6.7} Let $\Lambda$ be the lattice of highest weights of the ring $R_a$. Let $\Lambda_{even}$ (resp. $\Lambda_{odd}$) be the projection of the weight of $\GL(\FF)$ onto the weight of $\GL(\FF_{even})$ (resp. $\GL(\FF_{odd})$), where 
$$\GL(\FF)_{even} =\prod_{i\ even}\GL(F_i), \GL(\FF)_{odd} =\prod_{i\ odd}\GL(F_i).$$
Then the projections $\Lambda\rightarrow\Lambda_{even}$ (resp. $\Lambda\rightarrow\Lambda_{odd}$ are isomorphisms.
In other word every weight in $R_a$ is uniquely determined by its even and odd parts.
\endproclaim

Let us specialize to the case $n=3$. We will use slightly different notation.

 The incidence variety $Y_a$ giving a modification of variety $X_a:=Spec\ R_a$.
 
 $Y_a$ is a subset of $X_a\times Grass(r_3 ,F_2)\times Grass(r_2, F_1)\times Gras(r_1, F_0)$ consisting of tuples
    $ ((d_3,d_2,d_1,a_2,a_1), R_2, R_1, R_0)$ such that
    \item{}
   \item{1.} $\ (d_3, d_2, d_1, a_2, a_1)\in X_a$,
   \item{}
   \item{2.} $Im(a_i)\subset\bigwedge^{r_i}R_i, Im\ d_i\subset R_i\subset Ker\ d_{i-1}$ for $i=0,1,2$,
   \item{}
   \item{3.} For the induced maps $d'_i :Q_i:=F_i/R_i\rightarrow R_{i-1}$ and $a'_i\in\bigwedge^{r_i}R_i$, we have $det(d'_3)=a'_3$,
   $det(d'_2)=a'_3a'_2$, $det(d'_1)=a'_2a'_1$.
   \item{}

 One has natural projections $p_a :Y_a\rightarrow X_a$, $q_a:Y_a\rightarrow Grass$ where $Grass:=Grass(r_3 ,F_2)\times Grass(r_2, F_1)\times Gras(r_1, F_0)$.

 One gets the decomposition of $R_a$ (\cite{PW90}, section 1, formula (10) in \cite{W89})
 
 \proclaim{Proposition 6.8} We have
  $$R_a=\oplus_{a,b,c,\alpha,\beta, \gamma} S_{(a-b+c+\alpha_1 ,\ldots ,a-b+c+\alpha_{r_3-1},a-b+c)} F_3\otimes$$
  $$\otimes S_{(b-c+\beta_1 ,\ldots ,b-c+\beta_{r_2-1},b-c, -a+b-c, -a+b-c-\alpha_{r_3-1},\ldots ,-a+b-c-\alpha_1)}F_2\otimes$$
  $$\otimes S_{(c+\gamma_1 ,\ldots ,c+\gamma_{r_1-1}, c, c-b, c-b-\beta_{r_2-1},\ldots ,c-b-\beta_1)}F_1\otimes$$
  $$\otimes S_{(0^{r_0},-c, -c-\gamma_{r_1-1},\ldots ,-c-\gamma_1)}F_0.$$
  
  where we sum over all triples of natural numbers $a, b, c$ and triples of partitions $\alpha$, $\beta$, $\gamma$ such that $\alpha'_1<r_3$,
  $\beta'_1<r_2$ and $\gamma'_1<r_1$. 
   \endproclaim
   
   \proclaim{Corollary 6.9} The ring $R_a$ is a multiplicity free representation for the action of the group $\prod_{i=0}^3GL(F_i)$, so the variety $X_a$ is spherical.
   \endproclaim

In the case of $n=3$ we have over the ring $R_a$ $depth\ I(d_1)=1$, $depth\ I(d_2)=depth\ I(d_3)=2$.
In order to get acyclicity of $\FF_\bullet$ it is enough to raise the depth of $I(d_3)$ to $3$.
This can be done by killing the cycles in $\FF_\bullet$, in the Koszul complex on $I(d_3)$ and killing the higher direct image $R^1j_*(\OOO_{U_a})$ 
where $U_a=X_a\setminus V(I(d_3))$ and $j:U_a\rightarrow X_a$ is an inclusion. By the results of section 5 we know that vanishing of first homology groups of both complexes and of $R^1j_*(\OOO_{U_a})$ is equivalent.

We denote

$${\KKK}_a:\  0\rightarrow\bigwedge^0 {\KKK}\rightarrow\bigwedge^1{\KKK}\rightarrow \bigwedge^2 {\KKK}\rightarrow\bigwedge^3{\KKK}$$
  the beginning of the Koszul complex on $I(d_3)$, the ideal of maximal minors of $d_3$. Thus ${\KKK}:= \bigwedge^{r_3}F_3^*\otimes\bigwedge^{r_3}F_2\otimes_{\CC}R_a $. We treat ${\KKK}_a$ as a complex concentrated in degrees, $0$ to $3$ with differential of degree $-1$.

\proclaim{Proposition 6.10} We have the following isomorphisms.
\item{(1)} $H_1(\FF_a)= Ker (\bigwedge^3{\KKK}\otimes R^1j_*(\OOO_{U_a}){\buildrel{d_3\otimes 1}\over\rightarrow}\bigwedge^2{\KKK}\otimes R^1j_*(\OOO_{U_a}))$.
\item{(2)} $H_1(\KKK_a)= Ker (F_3\otimes R^1j_*(\OOO_{U_a}){\buildrel{d_3\otimes 1}\over\rightarrow}F_2\otimes R^1j_*(\OOO_{U_a}))$, i.e. it is the set of elements in $R^1j_*(\OOO_{U_a})$ annihilated by $I(d_3)$.
\endproclaim

Let us identify the generator of $H_1(\KKK_a)$.

\proclaim{Proposition 6.11} The module $H_1(\KKK_a)$ is generated by the image of the map
$$q_1:F_3^*\otimes\bigwedge^{r_1+1}F_1\rightarrow \bigwedge^2 (\bigwedge^{r_3}F_3^*\otimes\bigwedge^{r_3}F_2)\otimes R_a =$$
$$=S_{2^{r_3}}F_3^*\otimes\bigwedge^2(\bigwedge^{r_3}F_2)\otimes R_a .$$
There is the only one (up to nonzero scalar) nonzero equivariant map of this type.

In terms of the formula it is given by

$$h_t^*\otimes f_{i_1}\wedge\ldots\wedge f_{i_{r_1+1}}\mapsto\sum_{J, K} u_{I,J,K}g_J\otimes g_K$$
where the coefficient $u_{I,J,K}$ is given by the formula
$$\sum\pm \langle 1,\ldots ,{\hat t},\ldots ,r_3 |j_1\ldots ,{\hat j}_s,\ldots ,j_{r_3}\rangle_3\langle ( j_s, k_1,\ldots ,k_{r_3} )'| (i_1,\ldots ,i_{r_1+1})'\rangle_2 .$$
Here $I=(i_1,\ldots , i_{r_1+1})$, $J=(j_1,\ldots ,j_{r_3})$, $K=(k_1,\ldots ,k_{r_3})$, and $J'$ denotes the complement of the set $J$.
Also $\langle |\rangle_i$ denotes the minors of $d_i$ for $i=2,3$.
\endproclaim

\demo{Proof}

Let us look at possible equivariant maps $q_1$ as stated in the Proposition 6.11.
Looking at the weight corresponding to $\GL(F_1)$ we see that the only way such  map can occur is for the summand in $R_a$ having $b=c=0, \gamma=0$ and $\beta_1=\ldots =\beta_{r_2-1}=1$. Looking at the weight of $\GL(F_2)$
we need a trivial $SL(F_2)$ representation in $\bigwedge^2(\bigwedge^{r_3}F_2)\otimes S_{(1^{r_2-1},0^2,(-1)^{r_3-1})}F_2$.
It can happen only once, choosing the representation $S_{2^{r_3-1},1^2)}F_2$ in $\bigwedge^2(\bigwedge^{r_3}F_2)$.
Similar reasoning shows that the representation $F_3^*\otimes\bigwedge^{r_1+1}F_1$ cannot occur in
$\bigwedge^{r_3}F_3^*\otimes\bigwedge^{r_3}F_2\otimes R_a$.  Indeed, looking at the weight of $\GL(F_1)$ we see again that we need $b=c=\gamma=0$ and $\beta_1=\ldots =\beta_{r_2-1}=1$. Then looking at the weight of $\GL(F_2)$ we get a contradiction. 
Finally looking at the occurrence of  $F_3^*\otimes\bigwedge^{r_1+1}F_1$  in
$\bigwedge^3(\bigwedge^{r_3}F_3^*\otimes\bigwedge^{r_3}F_2)\otimes R_a$
we see it cannot happen, as the representation $S_{(3^{r_3-1}, 2,1)}F_2$ does not occur in 
$\bigwedge^3(\bigwedge^{r_3}F_2)$. 

\bigskip

There is another way to see that there is a cycle $q_1$ of the required form and that it generates $H_1({\KKK}_a)$.

We can calculate the higher direct image of $\OOO_{Z_a}$ with the map $a'_3$ inverted. This is done using Bott Theorem, applied to the quadruples of weights (corresponding respectively to $F_3, F_2, F_1, F_0$):

 $$((a-b+c+\alpha_1 ,\ldots ,a-b+c+\alpha_{r_3-1},a-b+c),$$
 $$ (b-c+\beta_1 ,\ldots ,b-c+\beta_{r_2-1},b-c, -a+b-c, -a+b-c-\alpha_{r_3-1},\ldots ,-a+b-c-\alpha_1),$$
 $$ (c+\gamma_1 ,\ldots ,c+\gamma_{r_1-1}, c, c-b, c-b-\beta_{r_2-1},\ldots ,c-b-\beta_1),$$
 $$ (0^{r_0},-c, -c-\gamma_{r_1-1},\ldots ,-c-\gamma_1))$$
  where we sum over all partitions $\alpha, \beta, \gamma$, $b,c\in\NN$ and $a\in\ZZ$, and then calculate the homology by Bott Theorem. We see that in $R^1j_*(\OOO_{U_a})$ we get the required representation for
  $$a=-2, b=c=\gamma =0, \beta_1=\ldots=\beta_{r_2-1}=1, \alpha_1=\ldots =\alpha_{r_3-1}=1.$$
  Moreover, if we increase $a$ by one to $-1$, there will be no corresponding representation in $R^1j_*(\OOO_{U_a})$,
  so our representation is annihilated by $I(d_3)$, so it gives an element in $H_1({\KKK}_a)$.
  
  The formula giving $q_1$ can be deduced from analyzing the way the equivariant map $q_1$ was constructed.
  Looking at the summand of $R_a$ we used it is clear it has to involve the products of $(r_3-1)\times (r_3-1)$ minors of $d_3$ and of $(r_2-1)\times (r_2-1)$ minors of $d_2$.\qed
  \enddemo
  
  \proclaim{Remark 6.12} Notice that we used the $\SL(F_2)\times \SL(F_1)$ equivariance instead of $\GL(F_2)\times \GL(F_1)$ equivariance to construct the map $q_1$. It is caused by the fact that under  the exact identification of weights in 
  the Lie algebra of type $T_{p,q,r}$ with the weights of $\GL(F_3)\times \GL(F_1)$ there is a copy of line bundle which centralizes ${\goth g}_0$ which acts in a nontrivial way.
  \endproclaim
  
  Let us compare the elements $q^{(1)}$ and $q_1$ as elements of $F_3\otimes R^1j_*(\OOO_{U_a})$ and $R^1j_*(\OOO_{U_a})$ respectively.
  Representation theory and Bott theorem show that they are related as follows.
  
  The map $q^{(1)}$ can be expressed as a composition
  $$\bigwedge^{r_1+1}F_1\otimes M_3^{-1}\otimes M_2\otimes M_1^{-1}{\buildrel{tr\otimes 1}\over\longrightarrow} F_3\otimes F_3^*\otimes \bigwedge^{r_1+1}F_1\otimes M_3^{-1}\otimes M_2\otimes M_1^{-1}\rightarrow$$
  $${\buildrel{q_1\otimes 1}\over\rightarrow} F_3\otimes R^1j_*(\OOO_{U_a}).$$
  
  Conversely, the map $q_1$ can be expressed in terms of $q^{(1)}$ as follows
  
  $$F_3^*\otimes \bigwedge^{r_1+1}F_1\otimes M_3^{-1}\otimes M_2\otimes M_1^{-1}{\buildrel{1\otimes q^{(1)}}\over\longrightarrow} F_3^*\otimes F_3\otimes R^1j_*(\OOO_{U_a}){\buildrel{ev\otimes 1}\over\rightarrow}  R^1j_*(\OOO_{U_a}).$$
  This follows from the fact that $R^1j_*(\OOO_{U_a})$  is multiplicity free.
  
  We conclude that adding to $R_a$ the entries of the cycle $b$ killing $q^{(1)}$ and cycle $p_1$ killing $q_1$, and performing the ideal transform with respect to $I(d_3)$ results in the same ring $R_1$. We denote $X_1= Spec(R_1)$ and $U_1=X_1\setminus V(I(d_3))$. Over the open set $U_a$ these rings are isomorphic to $\OOO_{U_a}$ with the variables corresponding to the first defect $F_3^*\otimes\bigwedge^{r_1+1}F_1$.
  
\head \S 7. The structure maps $p_i$. \endhead

\bigskip

 Let $\FF_\bullet$ be an acyclic complex of length three over a ring $R$. Let $\LL :=\LL(r_1+1, F_3, F_1)$ be the corresponding defect algebra.
 Finally, let
 $$0\rightarrow\bigwedge^0 {\KKK}\rightarrow\bigwedge^1{\KKK}\rightarrow \bigwedge^2 {\KKK}\rightarrow\bigwedge^3{\KKK}$$
 be the beginning of the Koszul complex on $I(d_3)$, the ideal of maximal minors of $d_3$. Thus ${\KKK}:= \bigwedge^{r_3}F_3^*\otimes\bigwedge^{r_3}F_2\otimes_{\CC}R $.
 
 In section 6 we constructed the map $p_1:\LL^*_1\rightarrow \bigwedge^1{\KKK}$ covering the cycle $q_1$. We continue to construct the higher maps $p_i$ ($i\ge 2$).
 
 \proclaim{Proposition 7.1} Let $\FF_\bullet$ be an acyclic complex of format $(f_0, f_1, f_2, f_3)$. There exists a structure map $p_2$
 making the following diagram commute.
 $$\matrix 0&\rightarrow&\bigwedge^0 {\KKK}&\rightarrow&\bigwedge^1{\KKK}&\rightarrow &\bigwedge^2 {\KKK}&\rightarrow&\bigwedge^3{\KKK}\\
 &&&&\ \ \ \uparrow p_{2}&&\ \ \ \uparrow \bigwedge^2 (p_1)&&\ \ \ \\\
 &&0&\rightarrow& \LL_{2}^*&\rightarrow& (\bigwedge^2\LL_1)^*&
 \endmatrix$$
 \endproclaim
 
 \demo{Proof} Since the upper row is an exact sequence, it is enough to check that the composition of the Koszul differential with $\bigwedge^2 (p_1)$ restricted to the image of $\LL_2^*$ is zero. This calculation is carried out in Theorem 2.9 from \cite{W89} (the map $p_1$ is denoted there by $b^\#$).
 \enddemo
 
 The defect of the map $p_2$ is equal to $\LL_2$.
 
 The definition of the defect algebra allows to introduce the higher maps $p_i$.
 
 \proclaim{Theorem 7.2}\cite{W89} Let $\FF_\bullet$ be an acyclic complex of format $(f_0, f_1, f_2, f_3)$.
 
 There exists a sequence of structure maps $p_i: \LL_i^*\rightarrow \bigwedge^1{\KKK}$ satisfying the following
 commutative diagram
 $$\matrix 0&\rightarrow&\bigwedge^0 {\KKK}&\rightarrow&\bigwedge^1{\KKK}&\rightarrow &\bigwedge^2 {\KKK}&\rightarrow&\bigwedge^3{\KKK}\\
 &&&&\ \ \ \uparrow p_{m+1}&&\ \ \ \uparrow q_{2,m+1}&&\ \ \ \uparrow q_{3,m+1}\\
 &&0&\rightarrow& \LL_{m+1}^*&\rightarrow& (\bigwedge^2\LL)^*_{m+1}&\rightarrow&(\bigwedge^3\LL)^*_{m+1}
 \endmatrix$$
 where $q_{2,m+1}=\sum (p_i\wedge p_j)$, $q_{3,m+1} =\sum (p_i\wedge p_j\wedge p_k)$.
 \endproclaim
 
 \demo{Proof} The upper row is an exact complex, and the diagram commutes by the definition of the Koszul differential and the maps in the lower row. The result follows by an elementary diagram chase.
 \enddemo

 The relation with the defect Lie algebra is that the defect (i.e. non-uniqueness) of each map $p_m$ is equal to $\LL_m$. Defect refers to the fact that $p_{m+1}$ is a lifting of certain cycles and we can modify $p_m$ by the map from $\LL_m^*$ to $R$,
 i.e. by an element of $\LL_m$. 
 
 The idea of the construction of the generic ring carried out in \cite{W89} is to build it up taking these symmetries into account. 
 More precisely, define $R_m$ to be the ring we obtain from $R_a$ 
 by adding generically the coefficients of the maps $p_1, \ldots ,p_m$, and then dividing by the appropriate relations (those that vanish when specializing to a splitting complex with arbitrary choice of the maps $p_1,\ldots ,p_m$, compare Lemma 2.4  \cite{W89}), and take the ideal transform with respect to $I(d_2)I(d_3)$.  We get the action of the Lie algebra $\LL /(\sum_{j>m}\LL_j)$ on such ring $R_m$ (Theorem 2.12, \cite{W89}).
 
 \proclaim{Definition 7.3} We define ${\hat R}_{gen}:=lim_{m} R_m$, $\FF_{gen} =\FF_a\otimes_{R_a} {\hat R}_{gen} .$
 \endproclaim
 
 Similarly, for every $m$ we have a diagram
 
  $$\matrix U_m:=Y_{m}\setminus p_{m}^{-1}(D_3 )&\buildrel{j'_{m}}\over\rightarrow&Y_{m}&\buildrel{q_{gen}}\over\rightarrow&Grass\\
   \downarrow\ p'_{m}&&\downarrow\ p_{m}&&\\
   X_{m}\setminus D_3&\buildrel{j_{m}}\over\rightarrow&X_{m}&&
   \endmatrix$$ 
   
   so finally, after including all $p_i$'s we get a  diagram
   
    $$\matrix U_{gen}:=Y_{gen}\setminus p_{gen}^{-1}(D_3 )&\buildrel{j'_{gen}}\over\rightarrow&Y_{gen}&\buildrel{q_{gen}}\over\rightarrow&Grass\\
   \downarrow\ p'_{gen}&&\downarrow\ p_{gen}&&\\
   X_{gen}\setminus D_3&\buildrel{j_{gen}}\over\rightarrow&X_{gen}&&
   \endmatrix$$ 
 
 Our goal is to show that $\RRR^1(j_{gen})_*\OOO_{X_{gen}\setminus D_3}=0$ proving that the complex $\FF_{gen}:=\FF_a\otimes_{R_a} \OOO_{X_{gen}}$ is the generic complex.
 
 \proclaim{Remark 7.4} Two observations will be useful in the future.
 \item{1)}The set $U_{gen}:=Y_{gen}\setminus p_{gen}^{-1}(D_3 )$ has a simple geometric interpretation.
 It is isomorphic to $U_0:=X_a\setminus p_a^{-1}(D_3)\times \oplus_{i>0}\LL_i$. Indeed, if the map $d_3$ splits, then each map $p_i$ splits into its defect and a map defined uniquely. Moreover, the affine space $\oplus_{i>0}\LL_i$ is clearly isomorphic to the open Schubert cell in the homogeneous space  $\GGG/\PPP$ where $\GGG$ is the Kac-Moody group associated to the graph $T_{p,q,r}$ and $\PPP$ is the parabolic associated to the simple root corresponding to the vertex $z_1$.
 \item{2)} The rings $R_m$ (and therefore the ring $R_{ten}$) are domains. Indeed, by construction the depth of the ideal $I(d_3)$ in these rings is $\ge 1$ and after inverting an $r_3\times r_3$ minor $D$ of $d_3$ we get a polynomial ring over $R_a[D^{-1}]$.
 \endproclaim

  \vskip 1cm

\head \S 8. The spectral sequence and the complexes $\KKK (\alpha,\beta ,s)_\bullet$ over $U(\LL)$. \endhead

   \bigskip
   
   The next step (section 3 in \cite{W89}) is the analysis of the spectral sequence allowing to calculate the cohomology of $\FFF_{gen}$.
   The spectral sequence is equivariant with respect to the group $\prod_{i=0}^3 GL(F_i)$. However the representations occurring
   in $U(\LL)$ do not contain the representations of $F_0$ and they contain only the maximal exterior power of $F_2$.
   Thus the $U(\LL)$-module structure will be preserved on the isotypic components of the group $\GG_{even}:=SL(F_2)\otimes GL(F_0)$.
   Let us also denote $\GG_{odd}:=SL(F_3)\otimes SL(F_1)$.
   
   Analyzing the isotypic components of the representations one reaches the following conclusion. The isotypic component of the cohomology of $\FFF_{gen}$
   is calculated as a cohomology of a complex $\KKK(\sigma ,\tau, t)$ of the form
   $$0\rightarrow K_0\rightarrow K_1\rightarrow K'_2\oplus K''_2$$
   where each term consists of a single irreducible representation of the group ${\hat\GG\LL}_{odd}:= GL(F_3)\times GL(F_1)$ tensored with $U(\LL)^*$.
   Dualizing we obtain the following statement.
   
   \proclaim{Theorem 8.1} (\cite{W89}, page 26, formula (38))
   All duals of isotypic components of the spectral sequence are the complexes of the form
   $$\matrix S_{(\sigma_1+t+u,\sigma_2,\ldots,\sigma_{r_3})}F_3\otimes \\
   \otimes S_{(\tau_1+t+u,\ldots, \tau_{r_1}+t+u,\tau_{r_1+1}+t, \tau_{r_1+2}+u,\tau_{r_1+3},\ldots ,\tau_{r_1+r_2})}F_1^*\otimes U(\LL)\\
   \oplus S_{(\sigma_1+t,\sigma_2+s,\sigma_3, \ldots,\sigma_{r_3})}F_3\otimes S_{(\tau_1+t+s,\ldots,\tau_{r_1+1}+t+s, \tau_{r_1+2},\ldots ,\tau_{r_1+r_2})}F_1^*\otimes U(\LL)\\
   \downarrow\\
   S_{(\sigma_1+t,\sigma_2,\ldots,\sigma_{r_3})}F_3\otimes S_{(\tau_1+t,\ldots,\tau_{r_1+1}+t, \tau_{r_1+2},\ldots ,\tau_{r_1+r_2})}F_1^*\otimes U(\LL)\\
   \downarrow\\
   S_{(\sigma_1,\sigma_2,\ldots,\sigma_{r_3})}F_3\otimes S_{(\tau_1,\ldots,\tau_{r_1+1}, \tau_{r_1+2},\ldots ,\tau_{r_1+r_2})}F_1^*\otimes U(\LL).
    \endmatrix$$
    Here the numbers $u$ and $s$ are uniquely determined by the triple $(\sigma,\tau, t)$ by equalities
    $$\sigma_2+s=\sigma_1+1, \tau_{r_1+2}+u=\tau_{r_1+1}+1.$$
   \endproclaim
   
   We denote the complex listed in the Theorem by $\KKK^* (\sigma,\tau, t)$.
   
   We have the following crucial consequence.
   
   \proclaim {Corollary 8.2} (\cite{W89}, Theorem 3.1)
   Assume that all the complexes $\KKK^* (\sigma,\tau, t)$ are exact at their middle term. Then $\RRR^1(j_{gen})_*\FFF_{gen}=0$ and therefore
   $$(j_{gen})_* \OOO_{X_{gen}\setminus D_3}=H^0 (Grass, \FFF_{gen})$$
   and the complex $\FF^{gen}_\bullet$ is acyclic over $(j_{gen})_* \OOO_{X_{gen}\setminus D_3}$, so it is the generic ring for our format.
   \endproclaim
   
    \vskip 1cm
    
    In the next section we will see that the complexes are indeed exact at the middle term by identifying them with the beginning part of certain parabolic $BGG$ 
    resolution.

\vskip 2cm

\head \S  9. Main result. \endhead

In this section we draw the consequences from previous cosiderations.
The main result is

\proclaim{Theorem 9.1} For every format $(f_0, f_1, f_2, f_3)$ there exists a generic pair
$$({\hat R}_{gen}, \FF^{gen}_\bullet):=((j_{gen})_* \OOO_{X_{gen}\setminus D_3}, \FF^a_\bullet\otimes_{R_a} (j_{gen})_* \OOO_{X_{gen}\setminus D_3}).$$
The generic ring ${\hat R}_{gen}$ is a general fibre of a flat family where a special fibre ${\hat R}_{spec}$ has a multiplicity free action of ${\goth g}(T_{p,q,r})\times {\goth gl}(F_2)\times {\goth gl}(F_0)$, where $f_3=r-1, f_2=q+r, f_1=p+q, r_1=p-1$.
If the algebra $\LL(r_1+1, F_3, F_1)$ is finite dimensional, then the generic ring ${\hat R}_{gen}$  is Noetherian.
\endproclaim

\demo {Proof} The complexes $\KKK^* (\sigma,\tau, t)$ from the Corollary 8.2. are identical to the part of $BGG$ complex identified in Proposition 4.4.
The partition $\alpha$ is just $\lambda$ restricted to the third arm of the graph. The partition $\beta$ is $\lambda$ restricted to the graph $A_{p+q-1}$ we get when we omit the third arm of the graph $T_{p,q,r}$. The number $t:=\lambda_{p+q}+1$. 
The differentials are the same because each component is nonzero and there is (up to a nonzero scalar) only one possible ${\goth gl}(F_3)\times {\goth gl}(F_1)$ map of free $U(\LL)$-modules in each case, so both differentials have to be the same.

Thus the complexes $\KKK^* (\sigma,\tau, t)$ are exact at the middle term so Corollary 8.2 assures that the complex $\FF^{gen}_\bullet$ is acyclic over ${\hat R}_{gen}$. 

Let us prove that the pair $({\hat R}_{gen}, \FF^{gen}_\bullet)$ has the universality property.
It was constructed by killing a series of cycles in the Koszul complex of $I(d_3)$. In every realization $(S, \GG_\bullet)$ where $S$ is Noetherian
and $\GG_\bullet$ is a resolution of format $(r_1 ,r_2, r_3)$ these cycles are boundaries. This and the universal property of the ideal transform give a homomorphism $\phi :{\hat R}_{gen}\rightarrow S$ such that
$\GG_\bullet = \FF_{gen}{}_\bullet\otimes_{{\hat R}_{gen}} S$.

This completes the proof that ${\hat R}_{gen}$ is indeed a generic ring.

\bigskip

To prove the part about the deformation, let us decompose the ring ${\hat R}_{gen}$ to the ${\goth gl}(F_2)\times {\goth gl}(F_0)$ isotypic components.

$${\hat R}_{gen}=\oplus_\mu {\hat R}_{gen,\mu}=\oplus_{\mu} S_{\phi(\mu)} F_0\otimes S_{\theta(\mu )}F_2\otimes V^*_{\lambda (\sigma(\mu), \tau(\mu), a)}$$

The cokernel of the complex 
$\KKK^* (\sigma,\tau, t)$ (i.e. the parabolic BGG complex)  is an irreducible highest weight module for ${\goth g}(T_{p,q,r})$, to the component ${\hat R}_{gen,\mu}$ acquires the structure of an irreducible ${\goth g}(T_{p,q,r})\times {\goth gl}(F_2)\times {\goth gl}(F_0)$-module. 
This action on $\oplus_\mu {\hat R}_{gen,\mu}$ is obviously multiplicity free. The problem is that this does not give the structure of the ${\goth g}(T_{p,q,r})\times {\goth gl}(F_2)\times {\goth gl}(F_0)$ module on the ring ${\hat R}_{gen}$ because the multiplication might not be ${\goth g}(T_{p,q,r})\times {\goth gl}(F_2)\times {\goth gl}(F_0)$ equivariant.

However for every two pieces ${\hat R}_{gen,\mu}$ and ${\hat R}_{gen,\nu}$ their product goes to the sum of several graded pieces with the extremal one being ${\hat R}_{gen,\mu +\nu}$.
We can deform the multiplication on ${\hat R}_{gen}$ by shrinking the other components of the product to zero. This gives us a new commutative algebra 
$${\hat R}_{spec}:=\oplus_{\mu\in\Lambda} {\hat R}_{gen,\mu}.$$

The connection between the rings ${\hat R}_{gen}$ and ${\hat R}_{spec}$ was explained in Grosshans lecture notes \cite{Gr97}, chapter 15. In theorem 15.14 Grosshans showed that there is an algebra $D$ which is a free $\CC[x]$ module such that
the general fibre of the resulting map 
$$\pi: Spec\ D\rightarrow \CC$$
over a point $z\in\CC$
 is isomorphic to $Spec\ {\hat R}_{gen}$ and the fibre over $0$ is isomorphic to ${\hat R}_{spec}$.
The next point is that the Cartan part of the multiplication map
$${\hat R}_{gen,\mu}\otimes {\hat R}_{gen,\nu}\rightarrow {\hat R}_{gen,\mu +\nu}$$
is not only ${\goth sl}(F_0)\times {\goth sl}(F_2)\times {\goth g}_+(T_{p,q,r})$-equivariant, but also
${\goth sl}(F_0)\times {\goth sl}(F_2)\times {\goth g}(T_{p,q,r})$-equivariant (so an epimorphism). The reason is as follows. It is well-known (see \cite{Ku02}, chapter X) that the homogeneous coordinate ring of the homogeneous space $\GGG/\PPP$ is a direct sum of irreducible representations 
$\oplus_{\lambda\in\Lambda}V(\lambda$ of irreducible representations $V(\lambda)$ ${\goth g}(T_{p,q,r}$ with $\Lambda$ consisting of all weights of type $\lambda(\sigma(\mu),\tau(\mu), a)$.  The multiplication map in this ring is just the Cartan multiplication $V(\lambda_1)\otimes V(\lambda_2)\rightarrow V(\lambda_1+\lambda_2)$, which is an epimorphism.

In order to compare the multiplications in the homogeneous coordinate ring of $\GGG/\PPP$ and the map ${\hat R}_{gen,\mu}\otimes {\hat R}_{gen,\nu}\rightarrow {\hat R}_{gen,\mu +\nu}$ we need one more fact.

 \vskip 1cm
    
     The complexes  $\KKK(\sigma ,\tau, t)$ we got in \cite{W89} as  isotypic components of the spectral sequence, before dualizing to get $\KKK^*(\sigma ,\tau, t)$ have another interpretation in terms of Grothendieck-Cousin complex introduced by George Kempf in \cite{Ke78}.
    The precise definitions of all the notions using in the remainder of this section can be found in \cite{Ku02}, chapter 9.
    
 The terms of the Grothendieck-Cousin complex are the local cohomology modules associated to the stratification of the homogeneous space $Z:={\GGG}/{\PPP}$ where ${\GGG}$ is the Kac-Moody group corresponding to $T_{p,q,r}$ and ${\PPP}$ is the parabolic subgroup corresponding to the simple root corresponding to the vertex $z_1$ (see section 4). The homogeneous space $Z$ has a stratification by Schubert cells and we denote by $Z_i$ the closed subset which is a union of all Schubert cells of codimension $\ge i$. Let ${\VVV}(\lambda(\sigma(\mu),\tau(\mu), a)$ be a homogeneous vector bundle on $Z$ corresponding to the weight $\lambda(\sigma(\mu),\tau(\mu), a)$.
    
    \proclaim {Proposition 9.2} (\cite{Ke78}, \cite{Ku02}, section 9.2)
    The isotypic component $\KKK(\sigma ,\tau, t)$ of the spectral sequence is the beginning part of the 
    Grothendieck-Cousin complex 
    $$0\rightarrow H^0_{Z_0/Z_1}(Z, {\VVV}(\lambda(\sigma(\mu),\tau(\mu), a))\rightarrow H^1_{Z_1/Z_2}(Z, {\VVV}(\lambda(\sigma(\mu),\tau(\mu), a))\rightarrow $$
    $$\rightarrow H^2_{Z_2/Z_3}(Z, {\VVV}(\lambda(\sigma(\mu),\tau(\mu), a)).$$
    \endproclaim
Now, looking at two weights $(\lambda (\sigma (\mu),\tau(\mu), a)$ and $(\lambda (\sigma (\nu),\tau(\nu), a')$ we see that bothe the Cartan part of their multiplication in ${\hat R}_{ten}$ and the multiplication of the elements of the kernels of Grothendieck-Cousin complexes from Proposition 9.2 are the same because they come from multiplication of sections on the open set $U_{ten}$ which, as we noted in Remark 7.4, is just the open Schubert cell, i.e. $Z_0\setminus Z_1$.

This allows us to prove  that if $T_{p,q,r}$ is a Dynkin diagram, then the rings ${\hat R}_{gen}$ and ${\hat R}_{spec}$ are Noetherian. 

We know that  the Lie algebra $\LL(r_1+1, F_3, F_1)$ is finite dimensional if and only if $T_{p,q,r}$ is a Dynkin diagram. In such case all irreducible highest weight modules for ${\goth g}(T_{p,q,r})$ are finite dimensional.
Therefore it is enough to show that the semigroup of weights occurring in ${\hat R}_{gen}$  is finitely generated.
But this semigroup is the semigroup of the terms in our spectral sequence which give the contribution to $H_0$. Thus we get the set of sextuples $(a,b,c,\alpha ,\beta,\gamma )$ with $a\in\ZZ$, $b, c\in \NN$ such that all the weights
$$(a-b+c+\alpha_1 ,\ldots ,a-b+c+\alpha_{r_3-1}, a-b+c),$$
$$ (b-c+\beta_1 ,\ldots, b-c+\beta_{r_2 -1}, b-c, -a+b-c, -a+b-c-\alpha_{r_3-1},\ldots, -a+b-c-\alpha_1) $$
$$(c+\gamma_1 ,\ldots ,c+\gamma_{r_1-1}, c, c-b, c-b-\beta_{r_2-1},\ldots ,c-b-\beta_1),$$
$$ (0^{f_0-r_1}, -c, -c-\gamma_{r_1-1},\ldots ,-c-\gamma_1) $$
are dominant. This  translates to the condition that $a\ge 0$, so our semigroup is finitely generated. 

\qqed
\enddemo

Let us summarize the properties of ${\hat R}_{spec}$.

\proclaim{Proposition 9.3}
We have an ${\goth sl}(F_0)\times {\goth sl}(F_2)\times {\goth g}(T_{p,q,r})$ decomposition
$${\hat R}_{spec}=\oplus_{\mu} S_{\phi(\mu)} F_0\otimes S_{\theta(\mu )}F_2\otimes V^*_{\lambda (\sigma(\mu), \tau(\mu), a)}$$
where $V_{\lambda}$ is the irreducible highest weight module of weight $\lambda$ for ${\goth g}(T_{p,q,r})$.
The dual $V^*_\lambda$ denotes the direct sum of duals of the weight spaces of $V_\lambda$. It is the highest weight representation for the opposite Borel subalgebra.
It is also irreducible.
The ring ${\hat R}_{spec}$ is a multiplicity free representation of ${\goth sl}(F_0)\times {\goth sl}(F_2)\times {\goth g}(T_{p,q,r})$. Its lattice of weights is saturated.
\endproclaim

\proclaim{Remark 9.4} I expect that the rings ${\hat R}_{gen}$ are also ${\goth gl}(F_0)\times {\goth gl}(F_2)\times {\goth g}(T_{p.q.r})$ equivariant.
In fact  it is enough to check the quadratic relations more precisely to see that they really hold in ${\hat R}_{gen}$. In every example analyzed below it is true.

\endproclaim

Let us exhibit the decomposition of ${\hat R}_{spec}$ explicitly.
For given $\alpha, \beta, t$ we define the weight $\lambda (\sigma, \tau, t)$ of ${\goth g}(T_{p,q,r})$ as follows.
We label the vertices of $T_{p,q,r}$ on the third arm by the coefficents of fundamental weights in $\sigma$, i.e.
$$\lambda_{p+q+i}=\sigma_{r-1-i}-\sigma_{r-i}$$
for $i=1,\ldots, r-2$.
We also label the vertices at the center and the first two arms by coefficients of fundamental weights in $\tau$,
i.e.
$$\lambda_0=\tau_p-\tau_{p+1},$$
$$ \lambda_i=\tau_{p-i}-\tau_{p-i+1}$$
for $1\le i\le p-1$, and
$$\lambda_i=\tau_i-\tau_{i+1}$$
for $i=p+1,\ldots, p+q-1$.
Finally, we put
$$\lambda_{p+q}=a.$$
We also set $t:=a+1$.

For a sextuple $\mu =(a,b,c,\alpha, \beta, \gamma )$ with $a\ge 0$, as in the decomposition of $(q_a)_*\FF_a$
we define 
$$\sigma (\mu):= (a-b+c+\alpha_1 ,\ldots ,a-b+c+\alpha_{r_3-1}, a-b+c),$$
$$\tau(\mu):=(c+\gamma_1 ,\ldots ,c+\gamma_{r_1-1}, c, c-b, c-b-\beta_{r_2-1},\ldots ,c-b-\beta_1),$$
$$\theta(\mu):= $$
$$=(b-c+\beta_1 ,\ldots, b-c+\beta_{r_2 -1}, b-c, -a+b-c, -a+b-c-\alpha_{r_3-1},\ldots, -a+b-c-\alpha_1),$$
$$\phi(\mu):= (0^{f_0-r_1}, -c, -c-\gamma_{r_1-1},\ldots ,-c-\gamma_1).$$

For the formats for which the algebra is not finite dimensional we do not have a Noetherian generic ring ${\hat R}_{gen}$. Still its multiplicity free structure  could be useful for applications.

\proclaim{Remark 9.5} In particular we proved the conjecture from \cite{PW90} stating that the generic ring ${\hat R}_{gen}$ constructed
there for the format $(f_0, f_1, f_2, f_3)=(1,n,n,1)$ is Noetherian. However the generators stated in \cite{PW90} in that case are not correct.
See the next section for more precise analysis of this case.
\endproclaim


\head  \S10. Properties of the rings ${\hat R}_{gen}$. \endhead

\bigskip

We start with the description of generators of the semigroup of weights of representations in ${\hat R}_{spec}$.
This will be useful when identifying the generators of ${\hat R}_{gen}$. This semigroup is the same as the semigroup of weights of $R_a$.
Therefore we have.

\proclaim{Proposition 10.1} The generators of the semigroup of weights in ${\hat R}_{spec}$ are as follows
\item{1)} $\alpha=(1^i), \beta=\gamma=a=b=c=0$, for $1\le i\le r_3-1$,
\item{2)} $a=1$, $\alpha=\beta=\gamma=b=c=0$,
\item{3)} $\beta=(1^j), \alpha=\gamma=a=b=c=0$, for $1\le j\le r_2-1$,
\item{4)} $b=1$, $\alpha=\beta=\gamma=a=c=0$,
\item{5)} $\gamma=(1^k), \alpha=\beta=a=b=c=0$, for $1\le k\le r_1-1$,
\item{6)} $c=1$, $\alpha=\beta=\gamma=a=b=0$.
\endproclaim

The multiplicity free property of the  ring ${\hat R}_{spec}$ also implies the following.

\proclaim{Theorem 10.2}
Assume that the ring ${\hat R}_{gen}$ constructed above is Noetherian (i.e. the graph $T_{p,q,r}$ is Dynkin).
Then $Spec\ {\hat R}_{gen}$ has rational singularities, in particular it is Cohen-Macaulay.
\endproclaim

\demo {Proof} As observed above the semigroup of weights of ${\hat R}_{spec}$  is saturated so the ring 
${\hat R}_{spec}$ is normal and therefore its spectrum is an affine spherical variety. This means that ${\hat R}_{spec}$
has rational singularities, by standard results on spherical varieties, see \cite{P12}, Corollary 4.3.15 and the following footnote.
Now ${\hat R}_{gen}$ also has rational singularities by the old result of Elkik \cite{E77}.
\qqed
\enddemo

Our next goal is to describe the (non-minimal) presentation of the ring ${\hat R}_{gen}$.
In order to do this let us recall the deformation technique of Grosshans \cite{Gr97}, chapter 15.

Given an irreducible affine  spherical variety $X$ for a reductive group $G$, with the coordinate ring 

$$K[X]=\oplus_{\lambda\in\Lambda} V_\lambda ,$$ 

where $\Lambda$ is some saturated lattice in the lattice of dominant integral weights for $G$, he constructed a flat deformation of $X$ with the special fibre $X_0$ where we also have

$$K[X_0]=\oplus_{\lambda\in\Lambda} V_\lambda$$

but the product in $K[X_0]$ is given by the Cartan product that multiplies two representations $V_\lambda$ and $V_\mu$ into their Cartan piece $V_{\lambda+\mu}$.

It is then known (Kostant's theorem, see \cite{Ku02}, section 10.1, especially Corollary 10.1.11) that the relations defining $K[X_0]$ are quadratic in the generators.

Applying this result we get

\proclaim{Theorem 10.3} The ring ${\hat R}_{gen}$ is generated by its subrepresentations ${\hat R}_{gen,\mu}$ corresponding to the generators $\mu$  of
the semigroup of weights for ${\hat R}_{spec}$  given by Proposition 10.1, with the relations that are quadratic in those generators.
This presentation is not minimal.
\endproclaim

The way to deal with relations is to use the following result from \cite{W89}.

\proclaim{Proposition 10.4}
A given $\prod_{i=0}^3 GL(F_i)$-equivariant set  $W$ of relations  holds in ${\hat R}_{gen}$ if and only if it vanishes when specializing to the splitting complex $\GG_\bullet$ of the format $(f_0, f_1, f_2, f_3)$ with the generic choice of the structure maps. 
\endproclaim

We conclude this section with the example of how one can work out the relations effectively.

\proclaim{Example 10.5} We analyze the presentation of ${\hat R}_{gen}$ for the format $(1, 4, 4, 1)$. The graph $T_{p,q,r}$ is a graph of type $D_4$
$$\matrix x_1&-&u&-&y_1\\ 
&&|&&\\
&&z_1&&
\endmatrix$$
We number the vertices as follows
$$\matrix 3&-&2&-&1\\ 
&&|&&\\
&&4&&
\endmatrix$$
The Lie algebra ${\goth g}(T_{p,q,r})$ is ${\goth so}(F_1\oplus F_1^*)$.

The weights of Proposition 9.3 are
$$\sigma(\mu)=(a-b+c),$$
$$\tau(\mu )=(c, c-b, c-b-\beta_2, c-b-\beta_1),$$
$$\theta(\mu )= (b-c+\beta_1, b-c+\beta_2, b-c, -a+b-c),$$
$$\phi (\mu)= (-c).$$
The weight $\lambda(\sigma(\mu),\tau(\mu), a)$ is

$$\matrix b&-&\beta_2&-&\beta_1-\beta_2\\ 
&&|&&\\
&&a&&
\endmatrix$$

The generators given by the Proposition 10.1. are. as follows (the generators of type 1) and 5) do not exist in this case).
\item{2)} $F_2^*\otimes V(\omega_4)$, 
\item{3)} $F_2\otimes V(\omega_1)$, $\bigwedge^2 F_2\otimes V(\omega_2)$,
\item{4)} $V(\omega_3)$,
\item{6)} $\CC$

Restricting to $F_1\oplus F_1^*$ we see that we get the generators
\item{2)} $F_2^*\otimes [\CC\oplus\bigwedge^2 F_1\oplus\bigwedge^4 F_1]=(d_3, p_1, v_2)$, 
\item{3)} $F_2\otimes V(\omega_1)=F_2\otimes (F_1^*\oplus F_1)=(d_2, c)$, $\bigwedge^2 F_2\otimes V(\omega_2)$ do not give new generators,
\item{4)} $V(\omega_3)=(F_1\oplus \bigwedge^3 F_1)=(a_2, p'_1)$,
\item{6)} $\CC=(a_1)$.

Here $d_3, d_2, d_1$ are differentials in $\FF_\bullet$, $a_2, a_1$ are the Buchsbaum-Eisenbud multipliers, $p_1:\bigwedge^2 F_1\rightarrow F_2$ is the second structure map, $c:F_1\otimes F_2\rightarrow F_3$ and $p'_1:\bigwedge^3 F_1\rightarrow F_3$ are gradede components of multiplication in $\FF_\bullet$, $v_2:\bigwedge^4 F_1\rightarrow  F_2$ is the structure map coming from the comparison of $\FF_\bullet$ and the Koszul complex on $a_2$ (see the next section).

We know that the relations are quadratic, so let us look at them in order. There are four basic representations (the representation of type 6) just adds the extra variable $a_1$). They are:
$$F_2^*\otimes V(\omega_4), F_2\otimes V(\omega_1), V(\omega_3)\ and\ \bigwedge^2 F_2\otimes V(\omega_2).$$
Let us look at the quadratic relations we get.

In degree $(2,0,0,0)$ we have
$$S_2(F_2^*\otimes V(\omega_4))=\bigwedge^2 F_2^*\otimes V(\omega_2)\oplus S_2 F_2^*\otimes [V(2\omega_4)\oplus\CC].$$
The first summand is the second representation of type 3). In the second summand the first part is a Cartan piece. The second part gives relation which is the factorization
 $d_3v_2=S_2(p_1)$.

The next case of degree  $(1,1,0,0)$ will be considered in some detail.
We look at the relations of type $F_2^*\otimes V(\omega_4)\otimes F_2\otimes V(\omega_1)$.
We have a tensor product decomposition $V(\omega_4)\otimes V(\omega_1)= V(\omega_1+\omega_4)\oplus V(\omega_3)$.
Of course we also have $F_2^*\otimes F_2 =adj(F_2)\oplus\CC$.
Out of four representations in the tensor product, two ($adj(F_2)\otimes V(\omega_1+\omega_4)$ and $\CC\otimes V(\omega_3)$) are in ${\hat R}_{gen}$ (as a Cartan piece and the generators of type 4)).

The remaining two representations form relations. Indeed we see from that representations of this type do not occur in ${\hat R}_{gen}$.

The representation $\CC\otimes V(\omega_1+\omega_4)$ consists of  relations that can be described as follows.
The representation $V(\omega_1)\otimes V(\omega_4)=(F_1^*+F_1)\otimes (\CC +\bigwedge^2 F_1+\bigwedge^4F_1)$, with the subrepresentation $V(\omega_3)=F_1+\bigwedge^3 F_1$.
The factor has four components: $F_1^*$ (giving relation $d_2d_3=0$), $S_{1,1,0,-1}F_1+S_{1,0,0,0}F_1$ (giving relations coming from the factorization $d_3 c=d_2p_1\pm x$), $S_{2,1,0,0}F_1 +S_{1,1,1,0}F_1$ (giving the relation between $d_2v_2$ and $cp_1$), and $S_{2,1,1,1}F_1$ (giving the relation $cv_2$).
 
The most interesting are the relations from the representation $adj(F_2)\otimes V(\omega_3)$. They are not readily seen from the relations of $R_a$.

We use Proposition 10.4.
Taking the splitting resolution in the form
$$d_3 =\bmatrix 0\\0\\0\\1\endbmatrix , d_2=\bmatrix 1&0&0&0\\0&1&0&0\\0&0&1&0\\0&0&0&0\endbmatrix , d_1=\bmatrix 0&0&0&1\endbmatrix $$
and denoting $\lbrace e_1, e_2, e_3, e_4\rbrace$ a basis in $F_1$, $\lbrace f_1, f_2, f_3, f_4\rbrace$ a basis in $F_2$, $g$ the basis element in $F_3$,
the generic multiplication is given by
$$e_1^.e_2=b_{12}f_4, e_1^.e_3=b_{13}f_4, e_2^.e_3= b_{23}f_4, $$
$$e_1^.e_4=-f_1+b_{14}f_4, e_2^.e_4= -f_2+ b_{24}f_4, e_3^.e_4=- f_3+b_{34}f_4$$
where $b_{i,j}$ are variables.
One calculates easily that
$$e_1^.f_1= 0, e_1^.f_2 = -b_{12}g, e_1^.f_3 = -b_{13}g, e_1^.f_4 = 0,$$
$$e_2^.f_1= b_{12}g, e_2^.f_2 = 0, e_2^.f_3 = -b_{23}g, e_2^.f_4 = 0,$$
$$e_3^.f_1= b_{13}g, e_3^.f_2 = b_{23}g, e_3^.f_3 = 0, e_3^.f_4 = 0,$$
$$e_4^.f_1= -b_{14}g, e_4^.f_2= -b_{24}g, e_4^.f_3= -b_{34}g, e_4^.f_4= g.$$
Moreover we have the multiplication
$e_1e_2e_3=0$, $e_1e_2e_4=b_{1,2}g$, $e_1e_3e_4=b_{1,3}g$, $e_2e_3e_4=b_{2,3}g$.
Finally we have the components of $v_2$ to be equal to $(v_2)_1=b_{2,3}$, $(v_2)_2=-b_{1,3}$, $(v_2)_3=b_{1,2}$, $(v_2)_4=b_{1,2}b_{3,4}-b_{1,3}b_{2,4}+b_{1,4}b_{2,3}$.

The  relations corresponding to $\CC \otimes V(\omega_1+\omega_4)$ involve the polynomials of the type
$$(v_2)_4 (d_2)_{1,1} = p_{23}^4c_4^1-p_{24}^4c_3^1+p_{34}^4 c_2^1$$
corresponding to the weight vector $f_4^*\otimes f_1\otimes e_2\wedge e_3\wedge e_4$ in $adj(F_2)\otimes\bigwedge^3 F_1\subset adj(F_2)\otimes V(\omega_3)$.

Using Proposition 10.4 we can prove these relations are satisfied. After specializing we get on both sides $b_{12}b_{34}-b_{13}b_{24}+b_{14}b_{23}$.

Relations of degree $(1,0,1,0)$ form the subrepresentation $F_2^*\otimes V(\omega_1)$ inside of $F_2^*\otimes V(\omega_4)\otimes V(\omega_3)$.
The first component $F_2^*\otimes\bigwedge^3 F_1$ contains  the commutativity relations of type $p_1a_2=p'_1d_3$.
The second component $F_2^*\otimes S_{2,1,1,1}F_1$ contains the equalities of type $p_1p'_1 =v_2a_2$.

Relations of degree $(1,0,0,1)$ are a subrepresentation in 
$$F_2^*\otimes V(\omega_4)\otimes\bigwedge^2F_2\otimes V(\omega_2)=$$
$$=[S_{1,0,0,0}F_2\oplus S_{1,1,0,-1}F_2]\otimes [V(\omega_2+\omega_4)\oplus V(\omega_1+\omega_3)\oplus V(\omega_4)].$$
Only $S_{1,1,0,-1}F_1\otimes V(\omega_2+\omega_4)$ survives in the ring ${\hat R}_{gen}$.
The relations come from the relations in degree $(1,1,0,0)$ by multiplying by $F_2\otimes V(\omega_1)$.
Indeed, we have that $V(\omega_1+\omega_4)\otimes V(\omega_1)$ maps epimorphically onto $V(\omega_4)\otimes V(\omega_2)$
and $V(\omega_3)\otimes V(\omega_1)$ maps onto $V(\omega_1+\omega_3)+V(\omega_4)$.

Relations of degree $(0,2,0,0)$ occur in
$$S_2(F_2\otimes V(\omega_1))= S_{2,0,0,0}F_2\otimes [V(2\omega_1)+\CC]\oplus S_{1,1,0,0}F_2\otimes V(\omega_2)$$
so only the representation $S_{2,0,0,0}F_2\otimes \CC$ is a relation. This is a relation of the form
$$\sum_{i=1}^4 (d_2)_{i,j}c_{k,i}+\sum_{i=1}^4 (d_2)_{i,k}c_{j,i}$$
 which can be checked by restricting to the splitting case.

Relations of degree $(0,1,1,0)$ occur in
$$F_2\otimes [V(\omega_1+\omega_3)\oplus V(\omega_4)]$$
and only $F_2\otimes V(\omega_4)$ give the relations.
Three components are $F_2\otimes \CC$ (giving the relation $a_2d_2=0$), $F_2\otimes\bigwedge^2F_1$ (giving the relations of type $d_2 p'_1+ca_2$), and $F_2\otimes\bigwedge^4 F_1$ (giving the relation of type $cp'_1$).

Relations of degree $(0,1,0,1)$ occur in
$$[S_{2,1,0,0}F_2\oplus S_{1,1,1,0}F_2]\otimes [V(\omega_1+\omega_2)\oplus V(\omega_3+\omega_4)\oplus V(\omega_1)].$$
The representations surviving in the ring are
$$S_{2,1,0,0}F_2\otimes V(\omega_1+\omega_2), S_{1,1,1,0}F_2\otimes V(\omega_3+\omega_4).$$
The relations come from the Cauchy formula (as $\bigwedge^3 V(\omega_1)=V(\omega_3+\omega_4)$ and $S_{2,1}V(\omega_1)= V(\omega_1+\omega_2)+V(\omega_1)$) and from multiplying the relation of degree $(0,2,0,0)$ by $F_2\otimes V(\omega_1)$. 

Relations of degree $(0,0,2,0)$ occur in
$$S_2 V(\omega_3)=V(2\omega_3)\oplus \CC$$
and only $\CC$ gives a relation which is of the type $a_2p'_1$, checked by restricting to the splitting case.

Relations of degree $(0,0,1,1)$ occur in
$$\bigwedge^2F_2\otimes [V(\omega_2+\omega_3)\oplus V(\omega_1+\omega_4)\oplus V(\omega_3)]$$
and two representations are in the relations. They come from multiplying the relations of degree $(0,1,1,0)$ by $F_2\otimes V(\omega_1)$.

Relations of degree $(0,0,0,2)$ occur in
$$[S_{2,2,0,0}F_2\oplus S_{1,1,1,1}F_2]\otimes [V(2\omega_2)\oplus V(2\omega_1)\oplus V(2\omega_3)\oplus V(2\omega_4)]\oplus$$
$$\oplus [S_{2,1,1,0}F_2\otimes [V(\omega_1+\omega_3+\omega_4)\oplus V(\omega_2)].$$
The surviving representations are $S_{2,2}F_2\otimes V(2\omega_2)$, $S_{2,1,1,0}F_2\otimes V(\omega_1+\omega_3+\omega_4)$ 
The relations involving $S_{2,2}F_2$ and $S_{2,1,1}F_2$ come from the relations involving Cauchy formula and the relations of degrees $(0,2,0,0)$ and $(0,1,0,1)$. The relations involving $S_{1,1,1,1}F_2$ come from the fact that 
$\bigwedge^4 V(\omega_1)=V(2\omega_3)+V(2\omega_4)$ and in the ring ${\hat R}_{gen}$ $S_{1,1,1,1}F_2$ occurs only with $V(\omega_3)$.

\vskip 1cm

We conclude this example with the comment on the action of ${\goth {gl}(F_2)}\times {\goth {gl}(F_0)}\times {\goth g}(T_{p,q,r})$ on ${\hat R}_{gen}$. The lowering operators for the action of ${\goth g}(T_{p,q,r})$ are the derivations acting on ${\hat R}_{gen}$, forming the algebra $\LL$. The subring of constants of these derivations is the ring $R_a$. In each representation
$$S_{\theta(\mu)}F_2\otimes S_{\phi(\mu)}F_0\otimes V(\lambda(\sigma(\mu),\tau(\mu), a)$$
 it picks the lowest component with respect to the induced grading on the module $V(\lambda(\sigma(\mu),\tau(\mu), a)$.
Let us consider the subring spanned by the highest graded components.
Recall that on the generators we had

\item{2)} $F_2^*\otimes [\CC\oplus\bigwedge^2 F_1\oplus\bigwedge^4 F_1]=(d_3, p_1, v_2)$, 
\item{3)} $F_2\otimes V(\omega_1)=F_2\otimes (F_1^*\oplus F_1)=(d_2, c)$, $\bigwedge^2 F_2\otimes V(\omega_2)$ do not give new generators,
\item{4)} $V(\omega_3)=(F_1\oplus \bigwedge^3 F_1)=(a_2, p'_1)$,
\item{6)} $\CC=(a_1)$.

So our subring contains $v_2, c, p'_1$ and $a_1$. Let us look at the sequence of maps:
$$\bigwedge^4 F_1\otimes {\hat R}_{gen}{\buildrel{v_2}\over\rightarrow}F_2\otimes {\hat R}_{gen}{\buildrel{c}\over\rightarrow}F_1^*\otimes {\hat R}_{gen}=\bigwedge^3 F_1\otimes{\hat R}_{gen}{\buildrel{p'_1}\over\rightarrow}F_3\otimes {\hat R}_{gen}.$$

The relations described above show this is a complex. Moreover, the structure map $p_1$ of this complex is our map $p_1$.
So the situation becomes self-dual. The resulting differential operators give the raising operators for the action of ${\goth g}(T_{p,q,r})$.
This is the simplest way to see that in this case ${\hat R}_{gen}$ indeed carries the action of ${\goth {gl}(F_2)}\times {\goth {gl}(F_0)}\times {\goth g}(T_{p,q,r})$ in this case.
\endproclaim

We conclude this section with some remarks on the expected generators of the ring ${\hat R}_{gen}$.
Before we start let us recall certain fact about the representations of ${\goth g}(T_{p,q,r})$.
Let us use the labeling on vertices by the labels $x_i, y_j, z_k, u$. We will use the convention that $u=x_p=y_q=z_r$.
For a given vertex $v$ we denote $V(\omega_v)$ the fundamental representation of ${\goth g}(T_{p,q,r})$
corresponding to the weight $\omega_v$ which has label $1$ at $v$ and label $0$ at all other vertices.

\proclaim {Proposition 10.6} For each $i$ ($1\le i\le p$) the representation $V(\omega_{x_{p-i}})$ occurs in $\bigwedge^iV(x_{p-1})$.
Similar statement is true for vertices $y_j$ and $z_k$.
\endproclaim

\demo {Proof} For the root systems of finite type this can be done by LiE computer algebra system. The general case follows from the fact that the highest weight
in $\bigwedge^i V(\omega_{x_{p-1}})$ is $\omega_{x_{p-i}}$.
\qqed
\enddemo

\proclaim{Remark 10.7} Let us look more closely at the generators of ${\hat R}_{gen}$ given in Proposition 10.1. 
\item{1)} The generators of type 1) for $i=1$ (or of type 2) if $r_3=1$) give is $F_2^*\otimes V(\omega_{z_{r-1}})$.
The first graded component of this representation (when considered as the representation of ${\goth g}_0(T_{p,q,r})$,
with the grading coming from the grading on ${\goth g}(T_{p,q,r})$), is  the map $d_3$. We denote the others by $p_i$.
They are indeed the maps $p_i$, all maps $p_i$ appear there. 
The later components are also tensors corresponding to maps $p_i$ occurring for the bigger formats (so for our format they sometimes do not have a defect).
This will be clear in the examples.
The maps of type $1)$ for higher $i$ should not give new generators: because of the Proposition 10.6 we expect these to be the minors of the generators of type $1)$ with $i=1$.
\item{2)} The generators of type 2): This is $\bigwedge^{r_2}F_2\otimes V(\omega_{z_1})$. These should also be the minors of generators of type $1)$ with $i=1$.
\item{3)} The generators of type 3) for $j=1$  give is $F_2\otimes V(\omega_{y_{q-1}})$.
The first graded component of this representation (when considered as the representation of ${\goth g}_0(T_{p,q,r})$,
with the grading coming from the grading on ${\goth g}(T_{p,q,r})$), is  the map $d_2$. We denote the others by $w_0, w_1,\ldots$.
In the most interesting case of $r_1=1$ the tensor $w_0$ is $F_2\otimes F_3^*\otimes F_2$ which gives part of the multiplicative structure on the resolution.
For higher $j$ the maps of type $2)$ should not give new generators: because of the Proposition 10.6 we expect these to be the minors of the generators of type $2)$ with $j=1$.
\item{4)} Here we  get the representation $V(\omega_{x_1})$. In the most interesting case $r_1=1$ of resolutions of the cyclic modules 
we get the first component to be $a_2$ and the second to be the tensor $F_3^*\otimes\bigwedge^3F_1$ which is the map $p'_1$ giving the other part of the multiplicative structure on the resolution. The other components are denoted $u_0, u_1,\ldots, $ with $u_0=p'_1$.
\item{5)} Again the generators for higher $k$ should just be the minors in the generators for $k=1$.
For $k=1$ we get a representation $F_0^*\otimes V(\omega_{x_{p-1}})$. these generators do not exist for $r_1=1$.
\item{6)} We get the representation $\bigwedge^{r_1}F_0^*\otimes V(\omega_{x_1})$. In the case $r_1=1$ this is just the map $a_1$. Otherwise this is the first component of our representation.
\endproclaim

\head  \S11. Examples \endhead

\bigskip

This section is devoted to presenting the explicit descriptions of the generic rings ${\hat R}_{gen}$ in the simplest cases.
We do this for the diagrams $T_{p,q,r}$ of types A and D.

\bigskip

\leftline{\bf 1. The graph $T_{p,q,r}$ is of type $A_n$} 

\bigskip

Formally there are three such cases: $p-1=0=r_1$, $q-1=0=r_2-2$, $r-1=0=r_3$. First and third cases are not legitimate,
so we consider the second. This means we have
$$dim\ F_3=r_3, dim\ F_2=r_3+2, dim\ F_1=r_1+2, dim\ F_0=r_1.$$

Before we go further let us make some remarks about the resolutions of this type.

Consider some free minimal acyclic complex
$$0\rightarrow G_3\buildrel{d_3}\over\rightarrow G_2\buildrel{d_2}\over\rightarrow G_1\buildrel{d_1}\over\rightarrow G_0$$
 over some Noetherian  ring
$S$, with 
$$dim\ G_3=r_3, dim\ G_2=r_3+2, dim\ G_1=r_1+2, dim\ G_0=r_1,$$
$$rank(d_3)=r_3, rank(d_2)=2, rank(d_1)=r_1.$$

There is a structure map $p_1:\bigwedge^{r_1+1}F_1\rightarrow F_2$ (the second structure map of Buchsbaum-Eisenbud, called $b$ in \cite{BE74}) and this map can be extended to the comparison map of complexes
$$\matrix 0&\rightarrow&G_3&\buildrel{d_3}\over\rightarrow&G_2&\buildrel{d_2}\over\rightarrow&G_1&\buildrel{d_1}\over\rightarrow&G_0\\
&&\uparrow\ p_1'&&\uparrow\ p_1&&\uparrow\ 1&&\uparrow\ 1\\
0&\rightarrow &\bigwedge^{r_1+2}G_1\otimes G_0^*&\buildrel{d_1^*}\over\rightarrow&\bigwedge^{r_1+1}G_1&\buildrel{q_1}\over\rightarrow&G_1&\buildrel{d_1}\over\rightarrow&G_0
\endmatrix$$

What happens is that after adding the coefficients of $p_1$ to the generic ring, the coefficients of $p_1'$ are added when taking the ideal transform with respect to $I(d_3)$. 

Now we come back to the generic ring we constructed.

The diagram $T_{p,q,r}$ looks like this

$$\matrix x_{p-1}&-&x_{p-2}&\ldots&x_{1}&-&u\\
&&&&&&|\\
&&&&&&z_{1}\\
&&&&&&|\\
&&&&&&\ldots\\
&&&&&&z_{r-2}\\
&&&&&&|\\
&&&&&&z_{r-1}
\endmatrix$$
with the distinguished root $z_1$.

The algebra ${\goth g}(T_{p,q,r})$ can be identified with ${\goth sl}(F_3\oplus \bigwedge^{r_1+1}F_1)={\goth sl}(F_3\oplus F_1^*)$.

Proposition 9.3 reduces to

$${\hat R}_{gen}=\oplus_{\mu} S_{\phi(\mu)} F_0\otimes S_{\theta(\mu )}F_2\otimes S_{\lambda (\sigma(\mu), \tau(\mu), a)} (F_3\oplus\bigwedge^{r_1+1}F_1)$$

where
for a sextuple $\mu =(a,b,c,\alpha, \beta, \gamma )$ with $a\ge 0$, as in the decomposition of $(q_a)_*\FF_a$
we define 
$$\sigma (\mu):= (a-b+c+\alpha_1 ,\ldots ,a-b+c+\alpha_{r_3-1}, a-b+c),$$
$$\tau(\mu):=(c+\gamma_1 ,\ldots ,c+\gamma_{r_1-1}, c, c-b, c-b-\beta_1),$$
$$\theta(\mu):= $$
$$=(b-c+\beta_1, b-c, -a+b-c, -a+b-c-\alpha_{r_3-1},\ldots, -a+b-c-\alpha_1),$$
$$\phi(\mu):= (0^{f_0-r_1}, -c, -c-\gamma_{r_1-1},\ldots ,-c-\gamma_1).$$

The weight $\lambda(\sigma(\mu), \tau(\mu), a)$ is the weight which corresponds to labeling the Dynkin diagram $T_{p,q,r}$ as follows
$$\matrix \tau_1-\tau_2&-&\tau_2-\tau_3&\ldots&\tau_{p-1}-\tau_p&-&\tau_p-\tau_{p+1}\\
&&&&&&|\\
&&&&&&a\\
&&&&&&|\\
&&&&&&\sigma_{r_3-1}-\sigma_{r_3}\\
&&&&&&\ldots\\
&&&&&&\sigma_2-\sigma_3\\
&&&&&&|\\
&&&&&&\sigma_1-\sigma_2
\endmatrix$$

The generators of the generic ring occur as the representations corresponding to the generators of the semigroup of highest weights of ${\hat R}_{gen}$. Some of them might be redundant, as the multiplication in
${\hat R}_{gen}$ is not a Cartan multiplication.

The representations giving the generators of the generic ring are as follows.
We number them according to Proposition 10.1.

\item {1)} with $i=1$,  i.e. $d_3$, i.e.  $F_3\otimes F_2^*$ corresponding to $\alpha$ being the first fundamental weight, $a=0$.
But this representation is really $F_2^*\otimes (F_3\oplus\bigwedge^{r_1+1}F_1)$, so it has coordinates $(d_3,  p_1)$.
If $\alpha$ is another fundamental representation, we get just the minors of the entries of the matrix $(d_3, p_1)$,
\item{2)} $a=1$, $\alpha =\beta =\gamma =0$, this corresponds to $r_3\times r_3$ minors of the matrix $(d_3, p_1)$.
\item{3)} with $j=1$, i.e.  $\beta_1=1$, $a=0$, $\alpha=\gamma=0$. This is the representation
$F_2\otimes F_1^*$ which gives the entries of $d_2$,
\item{4)} $b=1$, $a=c=\alpha=\beta=\gamma=0$, this representation is redundant, as it corresponds to $(r_3+2)\times (r_3+2)$ minors of $(d_3,p_1)$,
\item{5)} with $k=1$, i.e. $a=0$, $\alpha=\beta=0$, $\gamma$ is the first fundamental representation . We get the representation
$F_0^*\otimes (F_1\oplus F_3^*\otimes\bigwedge^3 F_1)$. The first part corresponds to $d_1$ and the second part  to the
structure map $p_1': F_0^*\rightarrow F_3$ constructed above. The higher fundamental representations $\gamma$ will just give the minors of the
two block matrix $(d_1, p'_1)$.
\item{6)} This is just the representation $a_1$.

We proved the following statement.

\proclaim{Theorem 11.1}  The generic ring ${\hat R}_{gen}$ for the resolution of the format $(r_3, r_3+2, r_1+2, r_1)$ is generated as an algebra over $R_a$
by the entries of  the structure maps $a_3, a_2, a_1$ and the structure maps $p_1, p_2$.
\endproclaim

Let us consider the relations between the generators.
The second structure theorem of Buchsbaum and Eisenbud (\cite{BE74}, section 6) says that we have the following commutative diagram

$$\matrix&&G_1^*=\bigwedge^{r_1+1}G_1&&\\
&\swarrow\ p_1&&\searrow\ a_1^*&\\
G_2=\bigwedge^{r_3+1}G_2^*&&\buildrel{d_2}\over\rightarrow&&G_1\\
&\searrow\ a_3^*&&\nearrow\ b^*&\\
&&G_2^*\endmatrix$$

This means that the generators $d_2$ and $a_1^*$ are redundant.
The defining relations for the generic ring are given by the two commutative diagrams and the relations coming from $R_a$.

A very interesting special case occurs when $r_1=1$.

\proclaim{Theorem 11.2}  The generic ring ${\hat R}_{gen}$ for the resolution of the format $(r_3, r_3+2, 3, 1)$ is a polynomial ring 
on the entries of the map $(d_3, p_1): F_3\oplus F_1^*\rightarrow F_2$ and the variable $a_1$.
\endproclaim

\demo{Proof} We know by the second structure theorem of Buchsbaum and Eisenbud that in this case the map $q_1=a_2$ is expressed through $d_3$ and $p_1$.
The only thing to show is to see that the map $p_1'$ is expressed through $d_3$ and $p_1$.
Looking at the diagram defining $p_1'$ is is not difficult to see that in the polynomial ring
$$Sym ((F_3\oplus F_1^*)\otimes F_2)\otimes Sym(F_0)$$
the entries of the map $p_1'$ will span the representation
$$\bigwedge^{r_3-1}F_3\otimes\bigwedge^{r_3+2}F_2^*\otimes \bigwedge^3 F_1^*\subset \bigwedge^{r_3-1}F_3\otimes\bigwedge^{r_3-1}F_2^*\otimes \bigwedge^3 F_1^*\otimes\bigwedge^3 F_2^*$$
in bidegree $(r_3-1, 3)$. So the generic ring ${\hat R}_{gen}$  is generated by the entries of $d_3$ and $p_1$.
We thus have
$${\hat R}_{gen}=Sym ((F_3\oplus F_1^*)\otimes F_2)\otimes Sym(F_0)$$
as one can see easily that there are no relations satisfied by $d_3$ and $p_1$ that vanish on any specialization to an acyclic complex.
Alternatively one could just prove using the exactness criterion of Buchsbaum and Eisenbud \cite{BE73} that the complex of length $3$
with differentials $d_3$, $d_2$, $d_1$ is acyclic over the polynomial ring. Also, one can identify the decomposition of ${\hat R}_{gen}$ to the irreducible representations with the one we get for
the polynomial ring using Cauchy formulas.

Just for completeness, let us give the formulas giving the differential $d_1$ in terms of $d_3$ $p_1$.
Denoting the matrix of $p_1 : F_1^*\rightarrow F_2$ by
$$B=\left(\matrix b_{1,1}&b_{1,2}&b_{1,3}\\
b_{2,1}&b_{2,2}&b_{2,3}\\
\ldots&\ldots&\ldots\\
b_{r_3+2,1}&b_{r_3+2,2}&b_{r_3+2,3}
\endmatrix\right).$$
Let us also denote
$$\Delta =\left( \matrix 0&\Delta_{1,2}&\Delta_{1,3}&\ldots&\Delta_{1,r_3+2}\\
-\Delta_{1,2}&0&\Delta_{2,3}&\ldots&\Delta_{2,r_3+2}\\
\ldots&\ldots&\ldots&\ldots&\ldots\\
-\Delta_{1,r_3+2}&-\Delta_{2, r_3+2}&-\Delta_{3, r_3+2}&\ldots&0\endmatrix\right)$$
the skew-symmetric matrix whose entry $\Delta_{i,j}$ is the maximal minor of the matrix of $d_3$ wth the $i$-th and $j$-th row omitted.
We have
$$d_1 =(a_1x_1, a_1x_2, a_1x_3)$$
where the entries $x_1, x_2, x_3$ are given by the matrix equation of skew-symmetric matrices
$$\left(\matrix 0&x_3&-x_2\\
-x_3&0&x_1\\
x_2&-x_1&0\endmatrix\right)= B^T\Delta B$$

\qqed
\enddemo

\proclaim{Remark 11.3} The case covered by Theorem 11.2 was already observed by Buchsbaum and Eisenbud in section 7 of \cite{BE74}. 
Notice that the present approach gives this result and the more general Theorem 11.1 from our main result 
just by analyzing the decomposition of ${\hat R}_{gen}$ into irreducible representations, without any additional considerations.
\endproclaim

\bigskip

\leftline{\bf 2. The graph $T_{p,q,r}$ is of type $D_n$} 
\bigskip

We have three possible cases.

\bigskip

\leftline{\bf a) $p=q=2$}

These are resolutions of the format $(1, 4, r_3+3, r_3)$.
The graph $T_{p,q,r}$ looks as follows

$$\matrix x_{1}&-&u&-&y_1\\
&&|&&\\
&&z_{1}&&\\
&&|&&\\
&&\ldots&&\\
&&z_{r-2}&&\\
&&|&&\\
&&z_{r-1}&&
\endmatrix$$
with the distinguished root $z_1$.
The corresponding Lie algebra is 
$${\goth g}(T_{p,q,r})={\goth so}(2(r_3+3)).$$ 
The orthogonal space in question is
$$U:= F_3\oplus \bigwedge^2 F_1\oplus F_3^*.$$

The generic ring is given by

$${\hat R}_{gen}=\oplus_{\mu} S_{\phi(\mu)} F_0\otimes S_{\theta(\mu )}F_2\otimes V^*_{\lambda (\sigma(\mu), \tau(\mu), a)} (F_3\oplus\bigwedge^2 F_1\oplus F_3^*)$$
where $V_\lambda$ denotes the irreducible representation of the special orthogonal Lie algebra of highest weight $\lambda$.

We have

$$\sigma(\mu)=(a-b+c+\alpha_1,\ldots,a-b+c+\alpha_{r_3-1},a-b+c),$$
$$\tau(\mu)=(c, c-b, c-b-\beta_2, c-b-\beta_1),$$
$$\theta(\mu)=(b-c+\beta_1, b-c+\beta_2, b-c, -a+b-c, -a+b-c-\alpha_{r_3-1},\ldots, -a+b-c-\alpha_1),$$
$$\phi(\mu)=(-c).$$
The weight $\lambda(\sigma(\mu),\tau(\mu),a)$ is given by labeling

$$\matrix b&-&\beta_2&-&\beta_1-\beta_2\\
&&|&&\\
&&a&&\\
&&|&&\\
&&\sigma_{r_3-1}-\sigma_{r_3}\\&&\\
&&|&&\\
&&\ldots&&\\
&&\sigma_2-\sigma_3&&\\
&&|&&\\
&&\sigma_1-\sigma_2&&
\endmatrix$$

We analyze the generators according to Proposition 10.1.

\item{1)} with $i=1$, $\alpha = (1,0^{r_3-2})$. We get the representation $F_2^*\otimes (F_3\oplus\bigwedge^2 F_1\oplus F_3^*\otimes\bigwedge^4 F_1)$.
The first component is just $d_3$, the second is $p_1$, the third one is $p_2$. It comes from the following factorization.
The usual comparison map

$$\matrix &0&\rightarrow&G_3&\buildrel{d_3}\over\rightarrow&G_2&\buildrel{d_2}\over\rightarrow&G_1&\buildrel{d_1}\over\rightarrow&R\\
&&&\uparrow\ p_1'&&\uparrow\ p_1&&\uparrow\ 1&&\uparrow\ 1\\
&R=\bigwedge^4 G_1&\buildrel{\delta}\over\rightarrow &\bigwedge^{3}G_1&\buildrel{\delta}\over\rightarrow&\bigwedge^{2}G_1&\buildrel{\delta}\over\rightarrow&G_1&\buildrel{\delta}\over\rightarrow&R
\endmatrix$$
gives a cycle $p_1'\delta$ which can be interpreted as a cycle

$$\matrix &0&\rightarrow&G_3&\buildrel{d_3}\over\rightarrow&G_2&\buildrel{d_2}\over\rightarrow&G_1&\buildrel{d_1}\over\rightarrow&R\\
&&&&&&&\uparrow\ u_2&&\\
&& &&&&&G_3^*&&
\endmatrix$$
which factorizes to $p_2:G_3^*\rightarrow G_2$ which is the third component of our representation. 
Alternatively, the third component is just the map $p_2$ defined in section 7.
\item{2)} These representations do not give new generators because all fundamental representations of ${\goth so}(2(r_3+3))$ up to the $r_3$-rd one are the exterior powers of $U$.
\item{3)} $F_2\otimes V^*(\omega_{r_3+2})(U)$ where $V(\omega_{r_3+2})(U)$ is half-spinor representation of ${\goth so}(U)$ which can be identified with
$\bigwedge^{even}(F_3)\otimes F_1\oplus \bigwedge^{odd}(F_3)\otimes F_1^*$.  The other half-spinor representation
$V(\omega_{r_3+3})(U)$ is half-spinor representation of ${\goth so}(U)$  can be identified with
$\bigwedge^{odd}(F_3)\otimes F_1\oplus \bigwedge^{even}(F_3)\otimes F_1^*$. 
We  note the following phenomenon which illustrates how tricky the identifications are.

Two half-spinor representations of the Lie algebra of type  $D_n$ are self-dual when $n$ is even and dual to each other when $n$ is odd.

This means in the generic ring we see  the representations
$$F_2\otimes (\bigwedge^{odd}(F_3)\otimes F_1\oplus \bigwedge^{even}(F_3)\otimes F_1^*)$$
when $r_3$ is even and
$$F_2\otimes (\bigwedge^{even}(F_3)\otimes F_1\oplus \bigwedge^{odd}(F_3)\otimes F_1^*)$$
when $r_3$ is odd.
Notice that in both cases we get in this representation the subrepresentation $F_2\otimes F_1^*$ which corresponds to $d_2$ (it has to be there after all). But notice also that in both cases we see the subrepresentation $F_2\otimes F_1\otimes F_3^*$ which corresponds to the multiplicative structure on our resolution. So we conclude that our ring indeed includes the tensor corresponding to the multiplicative structure, as it should.

One can rephrase it even better by saying that we will always see in ${\hat R}_{gen}$ the representations
$$F_2\otimes (\bigwedge^{odd}(F_3)^*\otimes F_1\oplus \bigwedge^{even}(F_3)^*\otimes F_1^*)$$

and

$$(\bigwedge^{even}(F_3)^*\otimes F_1\oplus \bigwedge^{odd}(F_3)^*\otimes F_1^*).$$

\item{4)}This gives another half-spinor representation $V^*(\omega_{r_3+3})(U)$ of ${\goth so}(U)$.
\item{5)} These representations do not exist as $r_1=1$,
\item{6)} This representation does not give new generators.

We conclude

\proclaim{Theorem 11.4} The generic ring ${\hat R}_{gen}$ for the format $(1, 4, r_3+3,r_3)$ is generated over $R_a$ by the entries of the representations of types 1), 3), 4). 
\endproclaim

\bigskip

\leftline{\bf b) $p=r=2$}

These are resolutions of the format $(1, r_2+1, r_2+1, 1)$. The graph $T_{p,q,r}$ looks like

$$\matrix x_{1}&-&u&-&y_1&-&\ldots&-&y_{q-1}\\
&&|&&&&\\
&&z_{1}&&&&
\endmatrix$$
with the distinguished root $z_1$. 

The Lie algebra
$${\goth g}(T_{p,q,r})={\goth so}(2(r_2+1)).$$
The orthogonal space in question is
$$U:= F_1\oplus F_1^*.$$

The formulas from section 9 reduce to

$${\hat R}_{gen}=\oplus_{\mu} S_{\phi(\mu)} F_0\otimes S_{\theta(\mu )}F_2\otimes V_{\lambda (\sigma(\mu), \tau(\mu), a)} (F_1\oplus F_1^*)$$
where $V_\lambda$ denotes the irreducible representation of the special orthogonal Lie algebra of highest weight $\lambda$.
Here we have
$$\sigma (\mu):= (a-b+c),$$
$$\tau(\mu):=(c, c-b, c-b-\beta_{r_2-1},\ldots ,c-b-\beta_1),$$
$$\theta(\mu):= (b-c+\beta_1 ,\ldots, b-c+\beta_{r_2 -1}, b-c, -a+b-c),$$
$$\phi(\mu):= (-c).$$

The weight $\lambda(\sigma(\mu), \tau(\mu ),a)$ corresponds to the labeling of the graph $T_{p,q,r}$ given as follows
$$\matrix \tau_{1}-\tau_{2}&-&\tau_{2}-\tau_{3}&-&\tau_{3}-\tau_4&-&\ldots&-&\tau_{q+1}-\tau_{q+2}\\
&&|&&&&\\
&&a&&&&
\endmatrix$$

We analyze the generators according to Proposition 10.1.

\item{1)} These representations do not exist as $r_3=1$,
\item{2)} This is a representation $F_2^*\otimes V(\omega_{n-1})$ where $V(\omega_{n-1})$ is half-spinor representation of ${\goth so}(F_1\oplus F_1^*)$. $V(\omega_{n-1})$ can be identified as $\oplus_{i\ge 0}\bigwedge^{2i} F_1$.
\item{3)} with $j=1$, this gives a multiplication $c: F_1\otimes F_2\rightarrow F_3$ on our resolution,
\item{4)}This gives another half-spinor representation $V(\omega_n)$ of ${\goth so}(F_1\oplus F_1^*)$ which can be identified as $\oplus_{i\ge 0}\bigwedge^{2i+1} F_1$.
\item{5)} These representations do not exist as $r_1=1$,
\item{6)} This representation does not give new generators.

Let us now describe how the representations from 2), 4) can be constructed from explicit cycles. This idea goes back to Bruns.

We construct a sequence of structure maps as follows.  
Let
$$0\rightarrow R{\buildrel{d_3}\over\rightarrow} G_2{\buildrel{d_2}\over\rightarrow} G_1{\buildrel{d_1}\over\rightarrow} R$$
be an acyclic complex over a Noetherian ring $R$. We have a comparison map from the Koszul complex

$$\matrix &0&\rightarrow&R&\buildrel{d_3}\over\rightarrow&G_2&\buildrel{d_2}\over\rightarrow&G_1&\buildrel{d_1}\over\rightarrow&R\\
&&&\uparrow\ p_1'&&\uparrow\ p_1&&\uparrow\ 1&&\uparrow\ 1\\
\ldots&\bigwedge^4 G_1&\buildrel{\delta}\over\rightarrow &\bigwedge^{3}G_1&\buildrel{\delta}\over\rightarrow&\bigwedge^{2}G_1&\buildrel{\delta}\over\rightarrow&G_1&\buildrel{\delta}\over\rightarrow&R
\endmatrix$$

Since $d_3 p_1'\delta =0$ we have $p_1'\delta=0$. This can be interpreted as a relation between the entries of $d_1$, more precisely, as a cycle
$u_2$ in the diagram, which can be then completed as follows.

$$\matrix &0&\rightarrow&R&\buildrel{d_3}\over\rightarrow&G_2&\buildrel{d_2}\over\rightarrow&G_1&\buildrel{d_1}\over\rightarrow&R\\
&&&\uparrow\ v_2'&&\uparrow\ v_2&&\uparrow\ 1&&\uparrow\ 1\\
\ldots&\bigwedge^6 G_1&\buildrel{\delta}\over\rightarrow &\bigwedge^{5}G_1&\buildrel{\delta}\over\rightarrow&\bigwedge^{4}G_1&\buildrel{u_2}\over\rightarrow& G_1&\buildrel{\delta}\over\rightarrow&R
\endmatrix$$

Since $d_3 v_2'\delta =0$ we have $v_2'\delta=0$. This can be interpreted as a relation between the entries of $d_1$, more precisely, as a cycle
$u_3$ in the diagram, which can be then completed as follows.

$$\matrix &0&\rightarrow&R&\buildrel{d_3}\over\rightarrow&G_2&\buildrel{d_2}\over\rightarrow&G_1&\buildrel{d_1}\over\rightarrow&R\\
&&&\uparrow\ v_3'&&\uparrow\ v_3&&\uparrow\ 1&&\uparrow\ 1\\
\ldots&\bigwedge^8 G_1&\buildrel{\delta}\over\rightarrow &\bigwedge^{7}G_1&\buildrel{\delta}\over\rightarrow&\bigwedge^{6}G_1&\buildrel{u_3}\over\rightarrow&G_1&\buildrel{\delta}\over\rightarrow&R
\endmatrix$$

In this way we construct the sequence of structure maps $v_i:\bigwedge^{2i}G_1\rightarrow G_2$ and $v_i':\bigwedge^{2i+1}G_1\rightarrow R$.

We define $v_1:=p_1$, $v_1':=p_1'$.

We have proved

\proclaim{Theorem 11.5} The generic ring ${\hat R}_{gen}$ for the format $(1, r_2+1, r_2+1,1)$ is generated over $R_a$ by the entries of the maps $c$, $v_i$, $v_i'$. 
\endproclaim

\proclaim{Remark 11.6} 
\item{a)} This case was analyzed in  \cite{PW90}, section 4. In fact we conjectured there (Conjecture 4.4)  that the generic ring for that format is Noetherian and we proposed the generators, (which did not include the $v_i'$'s) at the same time giving the elementary interpretation of $v_i$'s. In \cite{PW90} we worked assuming we automatically add 
to ${\hat R}_{gen}$ the fractions from ideal transforms with respect to $I(d_3)$. What happens is that adding the map $v_i$ automatically implies that in the ideal transform of the corresponding ring the entries of $v_i'$ appear.
\item{b)} The maps $v_i$ can be thought of as parts of higher maps $p_i$ existing for higher formats. Indeed, it is a pattern across all formats that generators called in Proposition 10.1 of type 2) (or of type 1) if rank of $F_3$ is equal to 1) involve the higher structure maps of type $p_i$.
\endproclaim

\leftline{\bf c) $q=r=2$}

These are resolutions of the format $(r_1, r_1+3,4,1)$.

$$\matrix x_{p-1}&-&x_{p-2}&\ldots&x_{1}&-&u&-&y_1\\
&&&&&&|&&\\
&&&&&&z_{1}&&\\
\endmatrix$$
with the distinguished root $z_1$.
The Lie algebra
$${\goth g}(T_{p,q,r})={\goth so}(2(r_1+3)).$$
The orthogonal space in question is
$$U:= F_1\oplus F_1^*.$$

The generic ring is given by

$${\hat R}_{gen}=\oplus_{\mu} S_{\phi(\mu)} F_0\otimes S_{\theta(\mu )}F_2\otimes V_{\lambda (\sigma(\mu), \tau(\mu), a)} (F_1\oplus F_1^*)$$
where $V_\lambda$ denotes the irreducible representation of the special orthogonal Lie algebra of highest weight $\lambda$.

$$\sigma (\mu):= (a-b+c),$$
$$\tau(\mu):=(c+\gamma_1 ,\ldots ,c+\gamma_{r_1-1}, c, c-b, c-b-\beta_2, c-b-\beta_1),$$
$$\theta(\mu):= $$
$$=(b-c+\beta_1, b-c+\beta_2, b-c, -a+b-c),$$
$$\phi(\mu):= (-c, -c-\gamma_{r_1-1},\ldots ,-c-\gamma_1).$$

The weight $\lambda (\sigma(\mu), \tau(\mu),a)$ is given by the labeling

$$\matrix \gamma_1-\gamma_2&-&\gamma_2-\gamma_3&\ldots&\gamma_{r_1-1}&-&b&-&\beta_2&-&\beta_1-\beta_2\\
&&&&&&&&|&&\\
&&&&&&&&a&&\\
\endmatrix$$

We analyze the generators according to Proposition 10.1.

\item{1)} These representations do not exist as $r_3=1$,
\item{2)} This is a representation $F_2^*\otimes V(\omega_{r_1+2})$ where $V(\omega_{r_1+2})$ is half-spinor representation of ${\goth so}(F_1\oplus F_1^*)$. $V(\omega_{r_3+2})$ can be identified as $\oplus_{i\ge 0}\bigwedge^{2i} F_1$.
\item{3)} with $j=1$, this gives a representation $F_2\otimes V(\omega_{r_3+3})$, where $V(\omega_{r_3+3})$ is another half-spinor representation
of ${\goth so}(U)$. It can be identified as $\oplus_{i\ge 0}\bigwedge^{2i+1} F_1$.
\item{4)}This gives the representation $\bigwedge^{r_1-3}(F_1\oplus F_1^*)$.
\item{5)} This, for $k=1$,  gives a representation $F_0^*\otimes (F_1\oplus F_1^*)$. Its first component is $d_1$.
\item{6)} This representation  does not give  new generators.

\bigskip

\head  \S12. Applications and symmetries \endhead

\bigskip

In this section we look at the possibilities of extending the structure theorems to perfect ideals of codimension three.
First we characterize the formats $(f_0, f_1, f_2, f_3)$ of Euler characteristic zero for which the minimal complexes
of a given format exist over a Noetherian local ring.
Throughout the section $(R, \goth m)$ will denote a local ring satisfying $depth(\goth m)\ge 3$.

We start with the perfect complexes resolving cyclic modules.
Thus we want to characterize the pairs $(n, l)$ such that there exist the acyclic complexes
$$0\rightarrow R^n\rightarrow R^{n+l}\rightarrow R^{l+1}\rightarrow R$$
resolving the cyclic modules $R/I$ with $depth(I)=3$.

\proclaim{Theorem 12.1}
Assume that one of the following conditions is satisfied
\item{a)} $l\ge 3, n\ge 2$,
\item{b)} $n=1$, $l$ even, $l>0$.

Then there exists a minimal  acyclic complex
$$0\rightarrow R^n\rightarrow R^{n+l}\rightarrow R^{l+1}\rightarrow R$$
resolving a cyclic module $R/I$ with $depth(I)=3$.
\endproclaim

\demo{Proof} We recall the basic construction involving linkage.
Let 
$$0\rightarrow R^n\rightarrow R^{n+l}\rightarrow R^{l+1}\rightarrow R$$
be an acyclic complex resolving a cyclic module $R/I$ with $depth(I)=3$.
Let $(x_1, x_2, x_3)$ be a regular sequence contained in $I$. We have the comparison map
from the Koszul complex on $(x_1, x_2, x_3)$ to our complex which gives rise to the diagram
$$\matrix 0&\rightarrow&R^n&\rightarrow&R^{n+l}&\rightarrow&R^{l+1}&\rightarrow&R\\
&&\uparrow\ \phi_3&&\uparrow\ \phi_2&&\uparrow\ \phi_1&&\uparrow\ 1\\
0&\rightarrow&R&\rightarrow&R^{3}&\rightarrow&R^{3}&\rightarrow&R.
\endmatrix$$

The dual of the mapping cone of the map $\phi$ then resolves $R/J$ where $J:=I:(x_1, x_2, x_3)$ is another ideal of depth $3$.

We can now make the following construction. Taking $x_1$ to be among the minimal generators of $I$, and $x_2, x_3$ to be in ${\goth m}^sI$ for such large $s$ that
$x_2$ and $x_3$ are not among the minimal generators of $I$ and the maps $\phi_2$ and $\phi_3$ have matrix entires in $\goth m$.
Then the dual of the mapping cone $C(\phi)_\bullet$ gives a minimal complex
$$0\rightarrow R^l\rightarrow R^{n+l+2}\rightarrow R^{n+3}\rightarrow R.$$
This means that the existence of our complex for the pair $(n,l)$ implies that similar complex exists for the pair $(l,n+2)$.

Now we perform the induction on $l+n$. Denoting the set of pairs $(n, l)$ satisfying the conditions $a)$, $b)$ by $\frak P$, we see that if
for a given $(n, l)$ the pair $(l-2,n)\in \frak P$ then by induction we establish the existence of the complex in question for the pair $(n, l)$.

Next we notice that if $n\ge 4, l\ge 4$ the pair $(l-2, n)$ is indeed in $\frak P$.
This means we need to deal with the remaining cases to provides the base for the induction.

We next observe that the case $l=3$, i.e. the case of an almost complete intersection is resolved in \cite{BE77}.

Next we deal with the case $n=3$.

We recall from \cite{BE77} that by using linkage argument and the structure theorem for Gorenstein ideals of codimension $3$ we can
construct for any $t\ge 2$ the minimal acyclic complex
$$0\rightarrow R^t\rightarrow R^{t+3}\rightarrow R^4\rightarrow R$$
resolving $R/J$ for a a perfect ideal $J$.

Taking a regular sequence $(x_1, x_2, x_3)$ in $J$ such that $x_1$ is among minimal generators of $J$ and $x_2,$ and $x_3$ are not, with the maps $\phi_2$, $\phi_3$ in the comparison map 
are minimal (this can be achieved by taking $x_2, x_3\in {\goth m}^s J$ for large $s$), we construct the minimal complex

$$0\rightarrow R^3\rightarrow R^{t+5}\rightarrow R^{t+3}\rightarrow R$$
satisfying the conditions of theorem for $n=3$.

Finally to deal with the case $n=2$ we again start with the resolution
$$0\rightarrow R^t\rightarrow R^{t+3}\rightarrow R^4\rightarrow R$$
resolving $R/J$ for a a perfect ideal $J$.
We choose a regular sequence $(x_1, x_2, x_3)$ in $J$ in such way that the first two elements are among the minimal generators of $J$, but the third is in the ${\goth m}^sJ$ for large $s$.

We can also arrange that the Koszul relation between $x_1, x_2$ is not a minimal generator of the syzygy module. This follows from the classification of the multiplicative structures on the resolutions of length 3
(see \cite{W89}).
Then the dual of the mapping cone of the comparison map gives a minimal complex of the type
$$0\rightarrow R^2\rightarrow R^{t+4}\rightarrow R^{t+3}\rightarrow R.$$
This takes care of the case $n=2$.
The case $n=1$ is the result of Buchsbaum and Eisenbud \cite {BE77}.
\qqed
\enddemo

\proclaim{Corollary 12.2}
Let $R$ be a local ring of depth $\ge 3$. Let $(f_0, f_1, f_2, f_3)$ be the format
of Euler characterisitc zero. Assume that we have $r_2>1$.
Then there exists a minimal acyclic free resolution of the format $(f_0, f_1, f_2, f_3)$
over $R$.
\endproclaim

\demo{Proof}  Let us handle the case $r_1=1$ first. By our assumption and Theorem 12.1
the required complex exists  unless $r_3=1$, $r_2$ is odd.
This means we need to construct the minimal acyclic complexes 
$$0\rightarrow R\rightarrow R^{2t}\rightarrow R^{2t}\rightarrow R.$$

We start with $R=\ZZ[X_1,\ldots ,X_{2t}]$ and we take the ideal $J$ generated by monomials
$J=(p_1,\ldots ,p_{2t})$ where
$$p_1 =X_1\ldots X_{2t-2}, p_2= X_2\ldots X_{2t-1},$$
$$p_i ={ {X_1\ldots X_{2t}}\over{X_{i-2}X_{i-1}}}$$
for $3\le i,\le 2t$.

One sees easily (using Buchsbaum-Eisenbud exactness criterion \cite{BE73}) that the free resolution of $R/J$ is given by the complex
$$0\rightarrow R(-2t)\rightarrow R^{2t}(-2t+1)\rightarrow R^{2t}(-2t+2)\rightarrow R$$
where the differentials are given by the matrices
$$d_3=\left(\matrix X_{2t}\\X_1\\X_2\\ \ldots\\ X_{2t-1}\endmatrix\right)$$

$$d_2=\left(\matrix X_{2t-1}&0&0&\ldots&0&-X_{2t}\\
-X_1&X_{2t}&0&\ldots&0&0\\
0&-X_2&X_1&\ldots&0&0\\
 \ldots\\ 0&0&0&\ldots&-X_{2t-1}&X_{2t-2}\endmatrix\right)$$
$$d_1 = (p_1, p_2, \ldots ,p_{2t}).$$

To specialize to an arbitrary ring $R$ we just need to find a sequence of elements $x_1,\ldots ,x_{2t}$ in $\goth m$ such that
each triple of elements $x_i, x_j, x_k$ forma a regular sequence, and specialize $X_i$ to $x_i$. The pecialized complex will still be acyclic by the exactness criterion \cite{BE73}.

To deal with the formats with $r_1>1$ we take the complex corresponding to the triple $(r_3, r_2, 1)$  and add to it $r_1-1$ copies of the complex  $0\rightarrow R{\buildrel x\over\rightarrow} R$
for some non-zerodivisor $x\in R$.

\proclaim{Remark 12.3} The condition $r_2>1$ is necessary by the first structure theorem of Buchsbaum and Eisenbud \cite{BE74}.
\endproclaim

\qqed
\enddemo

Next we explore the symmetries of our construction. We have in mind the possibility of generalizing the 
present approach to the structure theorems to perfect complexes and to see the effect of linkage on our construction.
 
Let us consider the format of the complex of length 3 with Euler characteristic zero.

Let us take $(f_0, f_1 , f_2 , f_3)=(p-1,p+q, q+r,r-1)$. Thus the sequence of ranks
is $(p-1, q+1, r-1)$.
Notice that this is the only format with Euler characteristic zero whose associated
defect Lie algebra $\LL$ is a parabolic corresponding to the graph $T_{p,q,r}$ with the 
distinguished root $z_1$.

For the dual complex the ranks will be $(r-1, q+1, p-1)$, so the corresponding algebra $\LL$ will be the parabolic associated to the Kac-Moody Lie algebra
corresponding to the same graph $T_{p,q,r}$ but with the distinguished root $x_1$.

Consider a perfect complex 
$$\GG_\bullet:0\rightarrow R^{p-1}\rightarrow R^{p+q}\rightarrow R^{q+r}\rightarrow R^{r-1}$$
 with ranks $(p-1, q+1, r-1)$ (i.e. such complex for which the dual complex is also acyclic).
We can costruct a perfect complex of another format    by using generalized linkage.
More precisely, consider the perfect complex of length $3$ of the format
$$\HH_\bullet :0\rightarrow R^{p-1}\rightarrow R^{p+1}\rightarrow R^{p+1}\rightarrow R^{p-1}$$
by choosing a $(p-1)\times (p+1)$ submatrix of the first differential $d_1$ of the complex $\GG_\bullet$ with the ideal of minors of grade three.
This can be always achieved by prime avoidance techniques.
There is an obvious map of complexes $\HH_\bullet\rightarrow \GG_\bullet$ and taking its mapping cone, dualizing and
minimalizing we get a complex of the format
$$0\rightarrow R^{q-1}\rightarrow R^{q+r}\rightarrow R^{p+r}\rightarrow R^{p-1}.$$
The  defect Lie algebra $\LL$ for the complexes of this format is again a parabolic for the Kac-Moody Lie algebra corresponding to the same graph $T_{p,q,r}$
but with the distinguished root $y_1$.

These symmetries could be viewed as a version of triality for the graphs $T_{p,q,r}$.
They give some hope that the approach will generalize to the perfect complexes.

We finish with several important remarks and questions:

\proclaim{Remarks 12.4}
\item{1.} The generic ring ${\hat R}_{gen}$ is a coordinate ring of some deformation of a homogeneous space for the group ${\goth g}(T_{p,q,r})\times GL(F_2)\times GL(F_0)$.
What is the geometric description of this object.
\item{2.} The types of the resolutions that could be related by linkage (in perfect case) come from the parabolic algebras in the same Kac-Moody Lie algebra.
Can we relate the perfect resolutions to Kac-Moody Lie algebras ? One possibility would be to use two defect Lie algebras: one related to $F_1$ and $F_3$ and the other (for the dual complex) built from $F_2^*$ and $F_0^*$. 
\item{3.} What is the relation of higher  terms of the spectral sequence with higher terms of the parabolic BGG complex ?
\item {4.} Can we get similar results for  modules of codimension 3 with a selfdual resolution (with an antisymmetric and symmetric middle matrix) ?
\item {5.} Can we get similar results for  Gorenstein ideals of codimension 4? They also are expected to have a manageable structure.
\item{6.} Does the approach generalize to resolutions of length bigger than 3 ? The depth raising procedure seems to work, we just need
more direct images $\RRR^ij_*$ to be zero. One possibility would be to have two Kac-Moody Lie algebras: one built from the even terms of $\FF_\bullet$ 
and the other from the odd terms of $\FF_\bullet$.
\endproclaim

\head  \S13. Acknowledgments \endhead

\bigskip

I would like to thank David Buchsbaum and David Eisenbud who introduced me to the subject of free resolutions and structure theorems, and who provided many insights.
I also had helpful discussions with Luchezar Avramov, Lars Christensen, Andrzej Daszkiewicz, Laurent Gruson, Craig Huneke, Witold Kra\'skiewicz, Andrew Kustin, Kyu-Hwan Lee, Matthew Miller, Alex Tchernev and Oana Veliche.

\bigskip\bigskip


\enddocument